\documentclass[12pt, a4paper]{article}

\usepackage{amsmath, amsthm, amssymb}  
\usepackage{fancyhdr}
\usepackage{graphicx}
\usepackage{enumerate}
\usepackage[top=2.5cm, bottom=2.5cm, left=2.5cm, right=2.5cm]{geometry}

\usepackage{mathrsfs} 
\usepackage{mathtools}
\usepackage{color}

\newtheorem*{definition}{Definition}
\newtheorem*{definition*}{Definition} 
\newtheorem{lemma}{\bf Lemma}
\newtheorem{proposition}{\bf Proposition}
\newtheorem{theorem}{\bf Theorem}
\newtheorem{corollary}{\bf Corollary}

\renewenvironment{proof}{\noindent {\bf Proof: }}{\rm\\}
\theoremstyle{definition}
\newtheorem{remark}{Remark}{\rm}
\newtheorem{example}[subsection]{Example}{\rm}

\usepackage{algorithm}
\usepackage{algorithmic}
\usepackage[T1]{fontenc}

\makeatletter
\newsavebox\myboxA
\newsavebox\myboxB
\newlength\mylenA

\newcommand*\xoverline[2][0.75]{%
    \sbox{\myboxA}{$\m@th#2$}%
    \setbox\myboxB\null
    \ht\myboxB=\ht\myboxA%
    \dp\myboxB=\dp\myboxA%
    \wd\myboxB=#1\wd\myboxA
    \sbox\myboxB{$\m@th\overline{\copy\myboxB}$}
    \setlength\mylenA{\the\wd\myboxA}
    \addtolength\mylenA{-\the\wd\myboxB}%
    \ifdim\wd\myboxB<\wd\myboxA%
       \rlap{\hskip 0.5\mylenA\usebox\myboxB}{\usebox\myboxA}%
    \else
        \hskip -0.5\mylenA\rlap{\usebox\myboxA}{\hskip 0.5\mylenA\usebox\myboxB}%
    \fi}
\makeatother

\newcommand{\overbar}[1]{\mkern 1.5mu\overline{\mkern-1.5mu#1\mkern-1.5mu}\mkern 1.5mu}

\usepackage{xparse}
\newcommand{\closure}[2][3]{%
{}\mkern#1mu\overline{\mkern-#1mu#2}}

\usepackage{auto-pst-pdf}
\usepackage{pst-plot}
\usepackage{pdftricks}
\begin{psinputs}
\usepackage{pstricks,pst-plot}
\end{psinputs}

\begin{document}

\title{Banach-Mazur game and open mapping theorem}
\author{Dominikus Noll
\thanks{Universit\'e de Toulouse, Institut de Math\'ematiques,  
118, route de Narbonne, 31062 Toulouse, France \newline
{\tt $\quad$ dominikus.noll$@$math.univ-toulouse.fr}}
}

\date{}

\maketitle

\begin{abstract}
We discuss a variant of the Banach-Mazur game which has applications to topological open mapping and closed graph theorems.

\vspace{.2cm}
\centerline{\bf Key Words}
\noindent
Open mapping theorem $\cdot$  Closed graph theorem $\cdot$ Banach-Mazur game $\cdot$ Baire space

\vspace{.2cm}
\centerline{\bf MSC 2020}

\centerline
{54E52 $\cdot$  54A10  $\cdot$  54C10 $\cdot$ 54B15 $\cdot$ 54D10}
\end{abstract}

\section{Introduction}
A Hausdorff topological space $E$ is called a {\it $B$-space}, respectively, a {\it $B_r$-space},  if every continuous, nearly open surjection, respectively bijection, $f$ from $E$ onto an arbitrary Hausdorff space $F$ is open. 
This definition honors V. Pt\'ak's famous open mapping theorem \cite{ptak_early,ptak,koethe} for linear operators between locally convex vector spaces, where 
he calls spaces $E$
satisfying the open mapping theorem {\it $B$-complete}, respectively,
{\it $B_r$-complete}. The concept has been extended to other categories, e.g.
Husain \cite{hu} calls a separated topological group $G$ a {\it $B$-group}  (a $B_r$-{\it group})  if every continuous nearly open surjective (bijective) homomorphism $f$ from $G$
into any separated topological group $H$ is open.  In the framework of  topological vector spaces, $B$- and $B_r$-completeness 
have been studied in \cite{adasch}.

 Weston  \cite{We}
proved that completely metrizable spaces are $B_r$-spaces, and Byczkowski and Pol \cite{BP} extended this to \v{C}ech-complete spaces. In
 \cite{open_mapping,baire_category}   
almost \v{C}ech complete spaces were shown to be $B_r$-spaces, while  locally compact spaces and 
Lindel\"off $P$-spaces are $B$-spaces
\cite{open_mapping}. Clearly every $B$-space is a $B_r$-space, while the converse is not true.

A natural question is whether, or in what sense, $B$- and $B_r$-spaces have to be {\it complete}. 
Pt\'ak's terminology points to the fact that locally convex spaces satisfying the open mapping theorem
are indeed {\it complete}.  But this is already different in the framework of topological groups, 
where  $B$- and $B_r$-groups may be incomplete and even of the first category. 
Since $H$-minimal topological spaces are $B_r$-spaces, and since Herrlich \cite{He} exhibits a first category $H$-minimal space, the topological case 
seems at first to resemble the situation in  groups.  However,  within metrizable spaces the situation is different, as all zero-dimensional
and all sub-orderable metrizable $B_r$-spaces are Baire  \cite{baire_category}. 
It is still an open question whether {\it all} metrizable $B_r$-spaces must be Baire. 

What is know is that 
there exist
metrizable $B_r$-spaces which do not contain any dense  completely metrizable subspace.  
Yet, all presently known constructions of metrizable $B_r$-spaces seem to require a somewhat strengthened form of Baire category.
This is in line with the fact
that
metrizability and Baire category  {\it alone} do  {\it not} suffice  to make a $B_r$-space. 

Presently we elucidate  this situation  a little further by discussing 
a novel variant of the Banach-Mazur game, where
$\beta$-defavorable spaces are $B_r$-spaces.  
This new variant is of interest in itself,  as the $\beta$-defavorable case implies Baire category not just of $E$, but of $E \times E$, while it still allows
$E\times E \times E$ to be of first category,  and  $E$ to be barely Baire. The $\alpha$-favorable case of the new game coincides  with the $\alpha$-favorable
case of the  classical Banach-Mazur game, but otherwise differences occur to the effect
that the gap between $\alpha$-favorability and $\beta$-favorability narrows. We  even identify classes of spaces where the game is determined.

It turns out that it is reasonable to consider, along with $B$- and $B_r$-spaces, a third intermediate concept.
We call a Hausdorff space $E$  a {\it $B_q$-space}, if every continuous nearly open surjection $f$ onto any Hausdorff space $F$ is a quotient map. 
In those categories where
quotient maps are by default open this offers nothing new, but $B$- and $B_q$-spaces
differ in the topological case. In fact, on closer regard, topological  $B$-spaces turn out a far more restrictive
class than one would expect from the situation in groups or topological vector spaces. Here
 $B_q$-spaces are a better behaved and sufficiently rich class to warrant further independent study.

The structure of the paper is as follows. After recalling the Banach-Mazur game in Section \ref{sectBM}, 
we introduce the novel variant in Section \ref{sectTandem}. Using the new game we obtain our first
open mapping theorem in Section \ref{sect_first}, followed by
a second open mapping theorem in Section \ref{sect_second} based on  the strong variant of the game.
The quotient mapping theorem is addressed in Section \ref{sect_quotient}. In Section \ref{sect_determined}
we resume the study of the new game and exhibit classes of spaces where it is determined.
Sections \ref{sect_barely},\ref{sect_barely_tau} on the other hand use the barely Baire spaces of \cite{fleissner} to show
that the new game is not determined in general. In Section \ref{sect_michael}
we have a glimpse at the Michael and  strong Choquet games and show that $\beta$-defavorability for those does 
not imply $\beta$-defavorability for our new game. In the final Section \ref{closed_graph} we prove a variant of the closed graph theorem,
where again two-person game characterizations play a crucial role.

\section{Preparations}
\label{prepare}
A mapping $f:E \to F$ is {\em nearly open} if for every $x\in E$ and every neighborhood $U$ of $x$ the set $\xoverline[.9]{f(U)}$
is a neighborhood of $f(x)$ in $F$.  Analogously, $f:E\to F$ is {\it nearly continuous} if for every $x\in E$ and every neighborhood $V$ of $f(x)$ in $F$, the set
$\xoverline[.9]{f^{-1}(V)}$ is a neighborhood of $x$.

Let $E$ be a Hausdorff space. A pair $(T,\phi)$ consisting of a tree $T = (T,\leqslant_T)$ of height $\omega$ and a mapping
$\phi$ from $T$ to the non-empty open subsets of $E$  is called a {\em web} on $E$ if 
\begin{itemize}
\item[$(w_1)$] $\{\phi(t): t \in T\}$ is a pseudo-base of $E$, i.e., every non-empty open $U \subset E$ contains some
$\phi(t)$.
\item[$(w_2)$] For every $t\in T$ the set $\{\phi(s): t <_T s\}$ is a pseudo-base of $\phi(t)$, i.e., every non-empty open
$V \subset \phi(t)$ contains some   $\phi(s)$ with $t <_T s$.
\end{itemize}

\begin{definition}[Completeness]
\label{def1}
{\rm 
A Hausdorff space $E$ is {\it $p$-complete} if it admits a web $(T,\phi)$ with the following property:
\begin{itemize}
\item[$(p)$]
 For every cofinal branch $b \subset T$
the intersection $\bigcap_{t\in b} \phi(t)\not=\emptyset$ is non-empty. 
\end{itemize}
%
The space $E$ is {\it $c$-complete} if it admits a web $(T,\phi)$ with the following property:
\begin{itemize}
\item[$(c)$]
 Every filter $\mathscr F$ containing all elements $\phi(t)$, $t\in b$,  of a cofinal branch $b\subset T$
has a cluster point in $\bigcap_{t\in b} \phi(t)$. 
\end{itemize}
}
\end{definition}

Pseudo-complete space in the sense of Oxtoby
\cite{ox} are $p$-complete, and $p$-completeness is also known as weak $\alpha$-favorability for a player with perfect information \cite{white}.
Almost \v{C}ech-complete space are $c$-complete, and coincide with the almost-complete spaces of \cite{almost_complete}.
Clearly $c$-completeness implies $p$-completeness.  Webs appear first in \cite{webbed,open_mapping}, and almost-sieves \cite{almost_complete}
are closely related. 
Our terminology generally follows \cite{engelking}, but we assume all spaces to be Hausdorff.

\section{Banach-Mazur game}
\label{sectBM}
The
Banach-Mazur game in a topological space $E$ is played by two players $\alpha$ and $\beta$ in the following way.
Player $\beta$ starts and chooses a non-empty open set $V_1$. Then player $\alpha$ chooses a non-empty open subset
$U_1 \subset V_1$. Next player $\beta$ chooses a nonempty open $V_2  \subset U_1$, and then $\alpha$ a non-empty open
$U_2 \subset V_2$, etc. Player $\alpha$ wins when $\bigcap_{n=1}^\infty U_n \not=\emptyset$, while player
$\beta$ wins when $\bigcap_{n=1}^\infty U_n=\emptyset$.

A strategy for player $\beta$ is a mapping, for simplicity also noted $\beta$, which for every sequence of non-empty open sets
$V_1,U_1,V_2,U_2,\dots,V_k,U_k$ of even length $2k$ chooses a non-empty open set $V_{k+1} = \beta(V_1,U_1,\dots,V_k,U_k) \subset U_k$.
This includes the sequence of length 0, where $\beta(\emptyset)=V_1 \not=\emptyset$.
Similarly, a strategy for player $\alpha$  is a mapping, noted $\alpha$,  which for every sequence
$V_1,U_1,\dots,V_k$ of non-empty open sets of odd length $2k-1$  chooses a non-empty open
$U_k = \alpha(V_1,U_1,\dots,V_k) \subset V_k$.
The play of $\alpha$ against $\beta$ is the sequence
$V_1\supset U_1 \supset V_2 \supset U_2 \dots$ satisfying
$V_1=\beta(\emptyset)$, $U_1=\alpha(V_1)$, $V_2 = \beta(V_1,U_1)$, $U_2=\alpha(V_1,U_1,V_2)$, etc. 
We call this the BM-game.

\begin{theorem}
\label{ox}
$E$ is a Baire space if and only if for every strategy $\beta$ in the {\rm BM}-game there exists a strategy $\alpha$
beating it. \hfill $\square$
\end{theorem}

The class of spaces $E$ where $\alpha$ has a winning strategy has been characterized  
by White \cite{white}, who called them weakly $\alpha$-favorable for a player with perfect information. These are
precisely the $p$-complete spaces  above. 

It is custom to call a space $\alpha$-favorable if player $\alpha$ has a winning strategy, and $\alpha$-defavorable if
it has not. The terms $\beta$-favorable and $\beta$-defavorable are used in the same sense. This terminology will also be used for variants of the
BM-game, and all variants will use perfect memory.


\section{Tandem Banach-Mazur game}
\label{sectTandem}

Now we introduce a variant of the Banach-Mazur game. We have two players $\alpha'$ and $\beta'$. 
Player $\beta'$ first chooses a non-empty open $V_1$, to which
player $\alpha'$ responds with a non-empty open $V_1' \subset V_1$. The second move of player $\beta'$  is 
a non-empty open $W_1$, to which player $\alpha'$ responds with a non-empty open $W_1' \subset W_1$. Then
player $\beta'$ switches back to the $V$-side and chooses a non-empty open $V_2 \subset V_1'$, to which $\alpha'$ responds with a non-empty open
$V_2' \subset V_2$. Then $\beta'$ chooses $W_2\subset W_1'$, and $\alpha'$ responds by $W_2' \subset W_2$, and so on. 
Players alternate between the $V$-side and the $W$-side, playing one move each on one side, before switching
to the other side and playing one move each there, etc.
The play therefore generates two nested sequences $V_1 \supset V_1' \supset V_2 \supset V_2' \supset \dots$
and $W_1 \supset W_1' \supset W_2 \supset W_2' \supset \dots$, but arranged in the following meandering order 

\begin{center}
\begin{tabular}{ccccccccc}
$\beta'$ && $\alpha'$ && $\beta'$ && $\alpha'$ && $\beta'$  \\
\hline
$V_1$ & $\to$ & $V_1'$               &                  & $V_2$        & $\to$               & $V_2'$  &                                    &$\dots$\\
            &                  &   $\downarrow$ &                  & $\uparrow$ &                                &       $\downarrow$      & &  $\uparrow$ \\
            &                  &    $W_1$            & $\to$ & $W_1'$         &                                &              $W_2$              &$\to$ & $W_2'$ \\
            \hline
           & & $\beta'$ && $\alpha'$ &&$\beta'$ && $\alpha'$ 
\end{tabular}
\end{center}
We
assume that both players have full information of the past
from both boards as time follows the arrows. 
Player $\alpha'$ wins when  both  $\bigcap_{n=1}^\infty V_n \not=\emptyset$ {\em and} $\bigcap_{n=1}^\infty W_n\not=\emptyset$, while
player $\beta'$ wins as soon as at least  one of these intersections is empty. We call this the {\it tandem Banach-Mazur game}, or for short,
the  ${\rm BM}'$-game.

\begin{proposition}
\label{prop1}
Player $\alpha$ has a winning strategy in the {\rm BM}-game if and only if player $\alpha'$ has a winning strategy
in the ${\rm BM}'$-game. 
\end{proposition}

\begin{proof}
1)
Suppose $\alpha$ is winning.
We construct a winning strategy $\alpha'$ by forgetting about the intertwined nature of the ${\rm BM}'$-game and by reacting to the moves
of $\beta'$ on the $V$- and $W$-board as if those were two non-cooperating players $\beta_V,\beta_W$ in the role of $\beta$ in the 
standard BM-game.
Since $\alpha$ wins against all strategies, it wins against these two, hence $\alpha'$ is winning as well.

2) Conversely, suppose $\alpha'$ is winning. We construct a winning strategy $\alpha$.  
Observe that
players $\beta'$ may choose  their moves such that $W_i \subset V_i'$, and also, $V_{i+1}\subset W_i'$, in which event the play in the ${\rm BM}'$-game 
will consist of  one single 
nested sequence and will correspond to a play in the BM-game. Therefore we may reinterpret $\beta$ as such a willful strategy $\beta'$.
The moves of the winning strategy $\alpha'$ on the $V$- and $W$-board may
then be translated back  to moves of a strategy $\alpha$ in the BM-game, which is winning since $\alpha'$ is.
\hfill $\square$
\end{proof}

Differences between the two games are expected when no winning strategy  $\alpha$ exists.
Let us observe that
in order to win player $\beta'$ has only to arrange for one of the nested sequences $V_i \supset V_{i+1}$ or
$W_i \supset W_{i+1}$ to have empty intersection. It is therefore clear that a winning strategy $\beta$ in the standard BM-game gives rise to a winning strategy
for $\beta'$ in the ${\rm BM}'$-game, by just concentrating on winning on one of the boards. 

\begin{proposition}
\label{beta1}
Suppose $\beta$ has a winning strategy in the {\rm BM}-game, then  $\beta'$ has a winning strategy in the ${\rm BM}'$-game. 
\end{proposition}

We will see later (Remarks \ref{easy} and  \ref{harder}) that the converse is not true, i.e., $\beta'$ may have a winning strategy, while $\beta$ has none.

An easy consequence of Proposition \ref{beta1}
is that if in $E$ the ${\rm BM}'$-game is $\beta'$-defavorable, 
then $E$ is a Baire space. But we have the following stronger

\begin{theorem}
\label{theorem3}
Suppose for every strategy $\beta'$ in the ${\rm BM}'$-game on $E$ there exists a strategy
$\alpha'$ which wins against $\beta'$. Then $E \times E$
is a Baire space.
\end{theorem}

\begin{proof}
We look at the Banach-Mazur game in the space $E \times E$ with players $\alpha,\beta$, and associate with them players
$\alpha',\beta'$ in the ${\rm BM}'$-game on $E$.

We may without loss of generality assume that $\beta$ plays with non-empty open boxes. Suppose $\beta(\emptyset)=V_1 \times W_1$.
We interpret $V_1$ as the first move of player $\beta'$ in the ${\rm BM}'$-game on $E$ on the $V$-board, to which player $\alpha'$ responds
by $V_1' \subset V_1$ nonempty. Normally, in the ${\rm BM}'$-game players $\beta'$ have now the possibility
to choose their move $W_1$ by taking $V_1,V_1'$  into account. But we take $W_1$ as the second component
of the first move of $\beta$ and let this be the second move of $\beta'$ in the ${\rm BM}'$-game, now
on the $W$-board. More formally,  with $p_1,p_2$ the projections on  first and second coordinate,
$\beta'(p_1(\beta(\emptyset)),V_1') = p_2(\beta(\emptyset))$ independently of the choice $V_1' \subset V_1 = p_1(\beta(\emptyset))$.
In other words, $\beta'$  wastes the option 
to take $V_1',V_1$ into account. 
Player $\alpha'$ reacts by choosing $W_1' \subset W_1$ based on all the previous information. 
Then we re-interpret $V_1' \times W_1'$ as the  move of player $\alpha$ in $E \times E$, i.e.,
$\alpha(V_1 \times W_1)=V_1' \times W_1' = \alpha'(V_1)\times \alpha'(V_1,\alpha'(V_1),W_1)$. We repeat this procedure in the following sweeps.

Once this wasteful strategy
$\beta'$  in the ${\rm BM}'$-game is defined, by hypothesis there exists $\alpha'$ beating $\beta'$.  The play so obtained
may now be read as a play $V_1\times W_1 \supset V_1' \times W_1' \supset V_2 \times W_2  \supset \dots$
between $\alpha$ and $\beta$ in the BM-game on $E \times E$. Since $\alpha'$ wins, we have $\bigcap_{i=1}^\infty V_i\not=\emptyset$ and
$\bigcap_{i=1}^\infty W_i\not=\emptyset$, hence of course $\bigcap_{i=1}^\infty V_i \times W_i \not=\emptyset$, hence $\alpha$ beats $\beta$. 
By Theorem \ref{ox}, $E\times E$ is a Baire space.
\hfill $\square$
\end{proof}

\begin{remark}
This construction generates a BM-strategy $\alpha$ on $E \times E$ winning against $\beta$, where  $\alpha$ plays
with open boxes in $E \times E$. It is not clear whether such a strategy exists when it is only known that $E\times E$ is Baire. We can
assume that $\beta$ plays with open boxes, but it is by no means clear whether this can be arranged for the $\alpha$ beating it.
\end{remark}

\begin{definition}[$\tau$-Baire space]
{\rm 
A topological space $E$ in which for every strategy $\beta'$ in the ${\rm BM}'$-game
there exists a strategy $\alpha'$ beating it is called
a {\it $\tau$-Baire space}.
}
\end{definition}

\begin{remark}
\label{easy}
It is well-known that there exist Baire spaces $E$ whose square $E \times E$ is no longer Baire (see \cite{fleissner} and Section \ref{sect_barely_tau}), and such a space
is Baire but not $\tau$-Baire. In Remark \ref{harder} we will see that even when $E \times E$ is Baire, this still does not mean that
$E$ is  $\tau$-Baire.
\end{remark}

\section{First open mapping theorem}
\label{sect_first}
A Hausdorff space $F$ is called a $\delta$-space if it admits  a web $(T,\phi)$ with 
the following property:
\begin{enumerate}
\item[$(d)$] For every cofinal branch $b \subset T$ the intersection $\bigcap_{t\in b} \phi(t)$ 
contains at most one point.
\end{enumerate}

Recall that a Hausdorff space is called {\it semi-regular}
if the family of regular-open sets is a basis for the topology;  \cite[p. 58]{engelking}. We are now ready to prove our first open mapping theorem:

\begin{theorem}
\label{theorem4}
Let $E$ be a
semi-regular $\tau$-Baire space, $F$ a $\delta$-space. If $f:E \to F$ is a continuous nearly open bijection,  then $f$ is open.
\end{theorem}

\begin{proof}
Let $(T,\phi)$ be a web on $F$ satisfying $(d)$.
We have to show that $f$ is open. 
Let $x\in E$ and $U$ a neighborhood of $x$. Using semi-regularity, choose an open neighborhood $V$ of $x$ with
$\xoverline[0.9]{V}^\circ \subset U$. It suffices to prove  $\xoverline[0.9]{f(V)}^\circ \subset f(U)$, as this will show openness at $x$.

Let $y\in \xoverline[0.9]{f(V)}^\circ$, $y = f(z)$. We will show $z\in \xoverline[0.9]{V}^\circ$, as this gives $z\in U$. By continuity of $f$
there exists an open $O$ with $z\in O$ and $f(O)\subset \xoverline[.9]{f(V)}^\circ$.
Proving $O \subset \xoverline{V}$ will now be sufficient. Take $w \in O$ and an arbitrary open neighborhood $W$ of $w$, which may be assumed 
to satisfy $W \subset O$. 
It remains to prove $V \cap W \not=\emptyset$, as this will show $w\in \xoverline[.9]{V}$.

We  define a strategy $\beta'$ in the ${\rm BM}'$-game on $E$. We start with the definition of $\beta'(\emptyset)$. 
Since $\{\phi(t): t\in T\}$ is a pseudo-base of $F$ by $(w_1)$, the set
$D_1=\bigcup\{\phi(t): t\in T, \phi(t) \subset \xoverline[.9]{f(V)}^\circ\}$ 
is open dense in $\xoverline[.9]{f(V)}^\circ$. On the other hand we have $f(W) \subset \xoverline{f(V)}^\circ$.
Now we apply Lemma \ref{lemma6} (below) with the choices
$G = D_1$, $H = \xoverline{f(V)}^\circ$, $O=V$, $U=W$. We conclude that $D_1 \cap \xoverline[.9]{f(W)}^\circ \cap f(V)\not=\emptyset$.
By the definition of $D_1$ we can  pick $t_1\in T$ with $\phi(t_1) \subset \xoverline[.9]{f(W)}^\circ$
and $x_1\in V$ satisfying $f(x_1) \in \phi(t_1)$. Now by continuity find an open $V_1$ with $x_1 \in V_1 \subset V$ and $f(V_1) \subset \phi(t_1)$.
Our move is $\beta'(\emptyset)=V_1$.

Let $V_1' \subset V_1$ be a potential response of player $\alpha'$. We have to define $\beta'(V_1,V_1')$. The set
$D_2= \bigcup\{ \phi(t): t\in T, t_1 <_Tt \}$ 
is open dense in $\phi(t_1)$ by $(w_2)$, while 
from
$V_1' \subset V_1$ we obtain $f(V_1') \subset \phi(t_1) \subset \xoverline[.9]{f(W)}^\circ$. 
We may therefore  apply Lemma \ref{lemma6} with the choices $G=D_2$, $H=\phi(t_1)$, $U=V_1'$ and $O=W$. The conclusion is that
$\xoverline[.9]{f(V_1')}^\circ$ intersects
$f(W) \cap D_2$. By the definition of $D_2$ choose $t_2\in T$ with $t_1 <_Tt_2$ and $\phi(t_2) \subset \xoverline[.9]{f(V_1')}^\circ$, together with $y_1\in W$ such that
$f(y_1) \in  \phi(t_2)$. By continuity find $W_1$ open with $y_1 \in W_1 \subset W$ satisfying $f(W_1) \subset \phi(t_2)$,
and let $\beta'(V_1,V_1') = W_1$ be our move.

Now let $W_1' \subset W_1$ be nonempty open, then we have to define
$\beta'(V_1,V_1',W_1,W_1')$. Note that
$D_3 := \bigcup \{\phi(t): t_2<_T t \}$ 
is dense in $\phi(t_2)$, while
$f(W_1') \subset \phi(t_2) \subset \xoverline[.9]{f(V_1')}^\circ$. This allows us to apply Lemma \ref{lemma6}
with the choices $H=\phi(t_2)$, $G=D_3$, $U=W_1'$ and $O=V_1'$. We conclude that
$\xoverline[.9]{f(W_1')}^\circ$ intersects $f(V_1') \cap D_3$. We can therefore pick $t_3\in T$ with $t_2 <_Tt_3$
satisfying $\phi(t_3)\subset \xoverline[.9]{f(W_1')}^\circ$ and $x_2 \in V_1'$ with $f(x_2) \in \phi(t_3)$. Pick an open $V_2$
with $x_2 \in V_2 \subset V_1'$ satisfying $f(V_2) \subset \phi(t_3)$, and put
$\beta'(V_1,V_1',W_1,W_1') = V_2$.

Let $V_2'\subset V_2$ be nonempty open. We have to define
$\beta'(V_1,V_1',W_1,W_1',V_2,V_2')$. The set
$D_4 = \bigcup \{ \phi(t): t_3<_Tt \}$ 
is dense in $\phi(t_3)$, while
$f(V_2')\subset f(V_2)\subset \phi(t_3)\subset \xoverline[.9]{f(W_1')}^\circ$. So we apply Lemma \ref{lemma6} with the choices 
$G=D_4$, $H=\phi(t_3)$, $O=W_1'$, $U=V_2'$. The consequence is that
$\xoverline[.9]{f(V_2')}^\circ$ intersects $f(W_1')\cap D_4$. We pick $t_4 \in T$, $t_3 <_Tt_4$ with $\phi(t_4) \subset \xoverline[.9]{f(V_2')}^\circ$
and $y_2 \in W_1'$ with $f(y_2) \in \phi(t_4)$. Then we choose an open set $W_2$ with $y_2 \in W_2 \subset W_1'$
such that $f(W_2)\subset \phi(t_4)$, and let
$\beta'(V_1,V_1',W_1,W_1',V_2,V_2') = W_2$ be our move.

Continuing in this way defines a strategy $\beta'$ in the ${\rm BM}'$-game, and since $E$ is $\tau$-Baire, we find a strategy $\alpha'$ winning
against $\beta'$. Let $V_1,V_1',W_1,W_1',V_2,V_2',W_2,W_2',\dots$ be their play. Then we find $v\in \bigcap_{i=1}^\infty V_i \not=\emptyset$
and $w\in \bigcap_{i=1}^\infty W_i\not=\emptyset$. In particular, $v\in V$ and $w\in W$. At the same time $f(v) \in f(V_i) \subset \phi(t_{2i-1})$
and $f(w) \in f(W_i) \subset \phi(t_{2i})$. Hence $f(v),f(w) \in \bigcap_{i=1}^\infty \phi(t_i)$, and since $t_1 <_T t_2 <_T t_3 <_T \dots$ is a cofinal branch in $T$,
we get $f(v) = f(w)$ from $(d)$. With $f$ being injective, this gives $v=w$, hence $V \cap W \not=\emptyset$. That ends the proof. 
\hfill $\square$
\end{proof}

\begin{lemma}
\label{lemma6}
Let $f:E \to F$ be continuous and nearly open. Let $O,U \subset E$ be open and $G,H \subset F$ open with $G \subset H \subset \xoverline[.9]{G}$.
Suppose $f(U) \subset H \subset \xoverline[.9]{f(O)}^\circ$. Then $G \cap \xoverline[.9]{f(U)}^\circ \cap f(O) \not=\emptyset$.
\end{lemma}

\begin{proof}
Observe that $G \cap \xoverline[.9]{f(U)}^\circ \not=\emptyset$. For had we $G \cap \xoverline[.9]{f(U)}^\circ =\emptyset$, then
also $\xoverline[.9]{G} \cap \xoverline[.9]{f(U)}^\circ=\emptyset$, hence $H \cap \xoverline[.9]{f(U)}^\circ = \emptyset$ due to $H \subset \xoverline[.9]{G}$.
But due to near openness that contradicts the assumption $f(U) \subset H$. Hence $G \cap \xoverline[.9]{f(U)}^\circ \not=\emptyset$.

From $f(U) \subset \xoverline[.9]{f(O)}^\circ$ follows $\xoverline[.9]{f(U)}^\circ \subset \xoverline[.9]{f(O)}^\circ$, and since $G \subset \xoverline[.9]{f(O)}^\circ$
anyway, the open set $G \cap \xoverline[.9]{f(U)}^\circ$ is contained in $\xoverline{f(O)}$. As it is nonempty by the above, it must
intersect the dense subset $f(O)$ of $\xoverline[.9]{f(O)}$.
\hfill $\square$
\end{proof}

%

Following Jayne and Rogers \cite{JR}, a topological space $F$ is {\it fragmentable}
if there exists a metric $d$ on $F$ such that for every $\epsilon > 0$ and every nonempty set $X \subset F$ there exists an open
set $U$ in $F$ such that $Y=X\cap U$ is nonempty and has diameter $\leq \epsilon$ with respect to $d$.
For a game-theoretic characterization of fragmentability see \cite{KM}. 

Following \v{C}oban et al. \cite{coban,coban2} a space $F$ is {\it open-fragmentable}
if the above condition applies to open sets $X\subset F$ only, in which case $Y$ is also open. The authors call this a {\it fos}-space
and obtain a related game-theoretic characterization, namely,
using plays $V_1 \supset W_1 \supset V_2 \supset W_2 \supset \dots$ where player $\alpha$ wins if
$\bigcap_{i=1}^\infty W_i$ consists of at most one point. Then $F$ is open-fragmentable iff player $\alpha$ has a winning strategy in this game
called the $FO$-game. See also \cite{wang} for further properties of this class.

\begin{lemma}
A space $F$ is open-fragmentable iff it is a $\delta$-space. 
\end{lemma}

\begin{proof}
1)
Suppose $F$ is open-fragmentable.
Define a tree $T$ of height $\omega$ as follows. The elements $t$ of $T$ are finite sequences of nonempty open sets
$t=(U_1,U_2,\dots,U_n)$ with $U_1\supset U_2 \supset \dots \supset U_n$ satisfying diam$(U_i) \leq 1/i$. The order relation is 
extension of sequences. The mapping $\phi$ is $(U_1,\dots,U_n) \to U_n$. It is clear that if $t_1 <_T t_2 <_T \dots$ is a cofinal branch, then this gives rise to a nested sequence
$U_1 \supset U_2 \supset \dots$ with $\bigcap_{i=1}^\infty U_i$ containing at most one point, hence property $(d)$ is guaranteed. 
We still have to prove that $(T,\phi)$ is a web.
Let us check property
$(w_2)$. Let $t=(U_1,\dots,U_n)\in T$. For every nonempty open $O \subset U_n$ there exists a nonempty open $U_{n+1} \subset O$
with diam$(U_{n+1})\leq 1/(n+1)$. Hence $\{U_{n+1}: (U_1,\dots,U_n,U_{n+1})\in T\}$ is a pseudo-base of the set $U_n=\phi(t)$. The argument for $(w_1)$ is similar.

2) For the converse, a web $(T,\phi)$ satisfying $(w_1),(w_2),(d)$ may obviously be used to define a winning strategy for $\alpha$ in the $FO$-game. Construction of a metric
which fragments open sets then follows as in \cite{coban2}.
\hfill $\square$
\end{proof}

\section{Second open mapping theorem}
\label{sect_second}
Now we consider the strong version of the ${\rm BM}'$-game, where player  $\alpha'$ wins the play
$(V_1,V_1',W_1,W_1',V_2,V_2',\dots)$ {\it strongly} against $\beta'$
if every filter $\mathscr F_V$ with $V_i\in \mathscr F_V$ for all $i$ has a cluster point
in $\bigcap_{i=1}^\infty V_i$, and every filter $\mathscr F_W$ with $W_i \in \mathscr F_W$ for all $i$ has
a cluster point in $\bigcap_{i=1}^\infty W_i$. 
If for every strategy $\beta'$ there exists a strategy $\alpha'$ winning strongly against $\beta'$, then we call $E$
a $\tau^*$-Baire space. The $c$-complete spaces of Section \ref{prepare} are those where players $\alpha$, or $\alpha'$, have a strong winning strategy,
so $c$-complete  spaces are $\tau^*$-Baire.

\begin{theorem}
\label{theorem5a}
Every semi-regular $\tau^*$-Baire space $E$ is a $B_r$-space.
\end{theorem}

\begin{proof}
Consider a continuous nearly open bijection $f:E \to F$.  We have to show that $f$ is open. 
Let $x\in E$ and $U$ a neighborhood of $x$. Using semi-regularity, choose an open neighborhood $V$ of $x$ with
$\xoverline[0.9]{V}^\circ \subset U$. It suffices to prove  $\xoverline[0.9]{f(V)}^\circ \subset f(U)$, as this will show openness at $x$.

Let $y\in \xoverline[0.9]{f(V)}^\circ$, $y = f(z)$. We will show $z\in \xoverline[0.9]{V}^\circ$, as this gives $z\in U$. By continuity of $f$
there exists an open $O$ with $z\in O$ and $f(O)\subset \xoverline[.9]{f(V)}^\circ$.
Proving $O \subset \xoverline{V}$ will now be sufficient. Take $w \in O$ and an arbitrary open neighborhood $W$ of $w$, which may be assumed    
to satisfy $W \subset O$. 
It remains to prove $V \cap W \not=\emptyset$, as this will show $z\in \xoverline[.9]{V}$.

We are going to define a strategy $\beta'$ in the ${\rm BM}'$-game on $E$. To start, observe that $\xoverline[.9]{f(W)}^\circ \cap f(V) \not=\emptyset$,
because $f(V)$ is dense in $\xoverline[.9]{f(V)}^\circ$, and $\xoverline[.9]{f(W)}^\circ$ intersects $\xoverline[.9]{f(V)}^\circ$, given that 
$f(W) \subset f(O)\subset \xoverline[.9]{f(V)}^\circ$.
Since $\xoverline[.9]{f(W)}^\circ \cap f(V)\not=\emptyset$,
we can find a nonempty open set $V_1$ with $f(V_1) \subset \xoverline[.9]{f(W)}^\circ \cap f(V)$.
Our first move is $\beta'(\emptyset)=V_1$.

Now  $\alpha'$ reacts to this by choosing $\emptyset\not= V_1' \subset V_1$, and we have to define $\beta'(V_1,V_1')$.
We have $\xoverline[.9]{f(V_1')}^\circ \cap f(W)\not=\emptyset$. Hence we can find
a nonempty open $W_1$ with
$W_1 \subset W$ and $f(W_1) \subset f(W)\cap \xoverline[.9]{f(V_1')}^\circ$. 
The corresponding move is now
$W_1=\beta'(V_1,V_1')$. 

Now player $\alpha'$ will react to this by providing $\emptyset \not= W_1' \subset W_1$. Since $f(W_1')\subset f(W_1)\subset  \xoverline[.9]{f(V_1')}^\circ$, we have
$\xoverline[.9]{f(W_1')}^\circ \cap f(V_1') \not=\emptyset$, so we pick a nonempty open  $V_2$ with
$V_2 \subset V_1'$ and $f(V_2) \subset
f(V_1') \cap  \xoverline[.9]{f(W_1')}^\circ$. Let $V_2=\beta'(V_1,V_1',W_1,W_1')$ be our move.

Next player $\alpha'$ chooses $\emptyset \not=V_2' \subset V_2$, and we have to define
$\beta'(V_1,V_1',W_1,W_1',V_2,V_2')$. From $f(V_2') \subset f(V_2)\subset \xoverline{f(W_1')}^\circ$
follows $\xoverline[.9]{f(V_2')}^\circ \cap f(W_1')\not=\emptyset$, so we can choose a nonempty 
open $W_2 \subset W_1'$ with $f(W_2) \subset \xoverline[.9]{f(V_2')}^\circ$.
The move is now $\beta'(V_1,V_1',W_1,W_1',V_2,V_2')=W_2$. Etc.

Having defined $\beta'$ in this way, let $\alpha'$ be a strategy which wins strongly against $\beta'$. Let
$V_1 \supset V_1' \supset V_2 \supset V_2' \supset \dots$ and $W_1 \supset W_1' \supset W_2 \supset W_2' \supset \dots$ be the two
nested sequences generated by the play. We have
\begin{itemize}
\item[i.] $f(W_i') \subset f(W_i) \subset \xoverline[.9]{f(V_i')}^\circ$, $W_i' \subset W_i \subset W_{i-1}'$, $W_1 \subset W$;
\item[ii.] $f(V_{i+1}') \subset f(V_{i+1}) \subset \xoverline[.9]{f(W_{i}')}^\circ$, $V_i' \subset V_i\subset V_{i-1}'$, $V_1 \subset V$.
\end{itemize}

Now choose a sequence $w_i\in W_i'$. 
As $\alpha'$ is strongly winning against $\beta'$, $w_i$ has a cluster point  $w\in \bigcap_{n=1}^\infty W_n\subset W$. 
Hence the sequence $f(w_i)$ has  cluster point $f(w)$. Now for every open neighborhood $G$ of $f(w)$ there exist
numbers $n(G,1) < n(G,2) < \dots$ such that $f(w_{n(G,i)}) \in G$ for all $i$.
By fact i. above $f(w_{n(G,i)}) \in \xoverline[.9]{f(V_{n(G,i)})}^\circ$, hence there exist $v_{n(G,i)}\in V_{n(G,i)}$ with
$f(v_{n(G,i)}) \in G$ for all $i$. Let  the set $\mathscr G$ of pairs $(G,i)$  be ordered by $(G,i) \preceq (G',i')$ iff $G' \subseteq G$ and $i' \geq i$.  Then
as $\alpha'$ is strongly winning the
net $\mathscr N =\langle v_{n(G,i)}: (G,i)\in \mathscr G\rangle$ has a cluster point $v\in \bigcap_{i=1}^\infty V_i \subset V$, hence $f(\mathscr N)$ has
cluster point $f(v)$. But by construction
the net $f(\mathscr N)$  converges to $f(w)$,  hence $f(w) = f(v)$. As $f$ is injective, we deduce $w=v$, and since $v\in V$, $w\in W$,
$V \cap W \not=\emptyset$ follows.
\hfill $\square$
\end{proof}

\begin{remark}
Consider a space $E$ possessing a web $(T,\phi)$ with the following property:
\begin{eqnarray*}
(\mu)&
\begin{array}{l}
\mbox{For every cofinal branch $b\subset T$, if $\bigcap_{t\in b} \phi(t) \not=\emptyset$ then every filter $\mathscr F$ containing all}\\ 
\mbox{$\phi(t)$, $t\in b$, has an 
cluster point in $\bigcap_{t\in b} \phi(t)$.}
\end{array}
\end{eqnarray*}
%
%
Then in $E$ the notions $\tau$-Baire and $\tau^*$-Baire coincide, and so do $p$- and $c$-completeness. 

In particular, the notions $\tau$-Baire
and $\tau^*$-Baire coincide in metrizable spaces, where $(T,\phi)$ with $(\mu)$ may be obtained by letting $T$ the tree of finite sequences
$(U_1,\dots,U_i)$ of non-empty open sets satisfying diam$(U_i) \leq 1/i$ and $\overline{U}_{i} \subset U_{i-1}$, ordered by extension of sequences, with $\phi$ denoting $(U_1,\dots,U_i) \to U_i$.
\end{remark}

\section{Inclusion and quotient mapping theorem}
\label{sect_quotient}
In this section we ask how to extend the open mapping theorem to more general continuous nearly open
mappings $f:E \to F$.
We consider the cases (a) $f$ injective and dense, i.e., $f(E)$ dense in $F$,    but no longer surjective, (b) $f$ surjective and no longer injective,
and (c) $f$ dense, and neither injective nor surjective.
We start with case (a), where the answer is easy.

\begin{proposition}
Let $E$ be a $B_r$ space and $f:E \to F$ a continuous, dense, injective and nearly open mapping into a Hausdorff space $F$.
Then $f$ is a homeomorphic embedding. 
\end{proposition}

\begin{proof}
It suffices to observe that $f$ considered as a mapping $f:E \to f(E)$ is still nearly open, as follows from Lemma \ref{lemma2} below. 
Then by the definition of a $B_r$-space
it is open, hence $E \simeq f(E)$.
\hfill $\square$
\end{proof}
 
\begin{lemma}
\label{lemma2}
Let $f:E\to F$ be a mapping such that $f(E)$ is dense in $F$. Then $f$ is nearly open if and only
if $f:E \to f(E)$ is nearly open.
\end{lemma}

\begin{proof}
1) Let $f$ be nearly open as a mapping $E \to f(E)$. Let $V \subset E$ be open. Then the closure of $f(V)$ in $f(E)$ is $cl_{f(E)} f(V) = \xoverline{f(V)} \cap f(E)$.
By assumption there exists  $O$ relatively open in $f(E)$ such that $f(V) \subset O \subset cl_{f(E)} f(V)$. Let $O = O_F \cap f(E)$ for $O_F$
open in $F$. Then $f(V) \subset O_F \cap f(E) \subset cl_{f(E)} f(V) \subset \xoverline{f(V)}$. Hence
$\xoverline[0.95]{O_F \cap f(E)} \subset \xoverline{f(V)}$. But $f(E)$ is dense in $F$, hence
$\xoverline[0.95]{O_F\cap f(E)} = \overline{O}_F$, proving $f(V) \subset O_F \subset \xoverline{f(V)}$. Hence $f:E \to F$ is nearly open.

2) Conversely, suppose $f:E \to F$ is nearly open. Let $V \subset E$ be open. Then $f(V) \subset O_F \subset \xoverline{f(V)}$ for
an open set $O_F$ in $F$. Therefore
$f(V) \subset O_F \cap f(E) \subset \xoverline{f(V)} \cap f(E)$. But $\xoverline{f(V)} \cap f(E) = cl_{f(E)} f(V)$, and $O=O_F \cap f(E)$
is relatively open in $f(E)$, hence $f:E \to f(E)$ is nearly open.
\hfill $\square$
\end{proof}

Concerning question (b) we have the following partial answer.

\begin{theorem}
\label{theorem5}
Let $E$ be $\tau^*$-Baire,  $f:E\to F$ a continuous and nearly open surjection.
Suppose $f$ is factorized as $f = g\circ h$ with $h:E \to G$ a continuous surjection onto a semi-regular space $G$ and $g:G\to F$ a continuous bijection. Then $g$ is open.
\end{theorem}

\begin{proof}
Note that $g$ is nearly open, because if $U \subset G$ is open, then
$h^{-1}(U)$ is open in $E$, hence $f(h^{-1}(U)) \subset \xoverline[.9]{f(h^{-1}(U))}^\circ$ by near openness of $f$. But clearly $f(h^{-1}(U))=g(U)$.

Let $x\in G$ and $U'$ a neighborhood of $x$. Using semi-regularity of $G$, choose an open neighborhood $V'$ of $x$ with
$\xoverline[0.9]{V'}^\circ \subset U'$. It suffices to prove  $\xoverline[0.9]{g(V')}^\circ \subset g(U')$, as this will show openness of $g$ at $x$.

Let $y\in \xoverline[0.9]{g(V')}^\circ$, $y = g(z)$. We will show $z\in \xoverline[0.9]{V'}^\circ$, as this gives $z\in U'$. By continuity of $g$
there exists an open $O'$ with $z\in O'$ and $g(O')\subset \xoverline[.9]{g(V')}^\circ$.
Proving $O' \subset \xoverline[.9]{V'}$ will now be sufficient. Take $w \in O'$ and an arbitrary open neighborhood $W'$ of $w$, which may be assumed    
to satisfy $W' \subset O'$. 
It remains to prove $V' \cap W' \not=\emptyset$, as this will show $z\in \xoverline[.9]{V'}$.

Now let $U=h^{-1}(U')$, $V=h^{-1}(V')$, $W=h^{-1}(W')$.  We define a strategy $\beta'$
on $E$ in exactly the same way as in the proof of Theorem \ref{theorem5a}. Following the proof all along will furnish
elements $v\in \bigcap_{i=1}^\infty V_i \subset V$ and $w\in \bigcap_{i=1}^\infty W_i \subset W$ for which $f(v)=f(w)$. Now a difference occurs,
as $f$ is no longer injective. But $g$ is, so from $f=g\circ h$ we obtain
$h(v)=h(w)$. Hence $h(v)=h(w) \in h(V)\cap h(W) =V' \cap W'$, and that was to be shown.
\hfill $\square$
\end{proof}

For the following recall that the {\it semi-regularization} of a space $E$, denoted $E_s$, is the point-set $E$ endowed with the coarser topology generated by the regular-open sets
of $E$
\cite[p. 58]{engelking}. The {\it regularization} of a space $E$, denoted $E_r$, is the point-set of $E$ endowed with the finest regular topology
coarser than the given one \cite{thomas}.  
Note that $E_r$ can  be defined explicitly using the ultra-closure operator of \cite{thomas}, where the author calls $E_r$  the associated regular space.

Suppose $E$ is regular, $\sim$ an equivalence relation on $E$,  $E\!/_{\!\sim}$ the quotient space, $\phi:E \to E\!/_{\!\sim}$
the quotient map. Since $E\!/_{\!\sim}$ is not necessarily regular, we use its regularization $(E\!/_{\!\sim})_r$ and consider $\phi$ as a mapping
$\phi:E \to 
(E\!/_{\!\sim})_r$. We call $\phi$ a {\it regular-quotient} map and $(E\!/_{\!\sim})_r$ the {\it regular-quotient}, 
because it preserves the following universal property of quotients:  if $g: (E\!/_{\!\sim})_r \to F$
is any mapping into a regular space $F$, then $g$ is continuous iff $f=g \circ \phi$ is continuous. In the same vein, we call a continuous surjection $f$
from a regular space $E$ onto a regular space $F$ {\it regular-quotient} if $F \simeq (E\!/_{\!\sim})_r$ for the equivalence relation 
$x \sim y $ iff $f(x)=f(y)$.

\begin{corollary}
\label{cor1}
Let $E$ be a regular $\tau^*$-Baire space, $F$ a regular space, and $f:E\to F$ a continuous nearly open surjection. Then
$f$ is a regular-quotient mapping, i.e. the topology on $F$ is the regularization, and also the semi-regularization, of the quotient topology.
\end{corollary}

\begin{proof}
Consider the equivalence relation $x\sim y$ iff $f(x)=f(y)$ and let $E\!/_{\!\sim}$ be the quotient space with the usual quotient topology, 
$\phi:E \to E\!/_{\!\sim}$ the quotient map,
$\widetilde{f}:E\!/_{\!\sim} \to F$ the continuous bijection satisfying $\widetilde{f} \circ \phi = f$. Then $E\!/_{\!\sim}$ is Hausdorff because $F$ is, but
$E\!/_{\!\sim}$ need not be semi-regular.
Let $\left(E\!/_{\!\sim}\right)_s$ be the semi-regularization of $E\!/_{\!\sim}$. Then
$\left(E\!/_{\!\sim}\right)_s$ is still Hausdorff and 
$\phi:E \to\left(E\!/_{\!\sim}\right)_s$ is continuous.  The point is now that $\widetilde{f}:\left(E\!/_{\!\sim}\right)_s \to F$ remains continuous due to regularity of $F$, 
see \cite[I.3.(3)]{katetov} or \cite[Prop. 2.2g]{porter}. 
We may therefore apply Theorem \ref{theorem5} with $G=\left(E\!/_{\!\sim}\right)_s$, $h=\phi$, $g=\widetilde{f}$,
which shows that $\widetilde{f}$ is a homeomorphism. That gives $\left(E\!/_{\!\sim}\right)_s \simeq F$, so $\left(E\!/_{\!\sim}\right)_s$ is regular.

Now by definition the regularization of the quotient $\left(E\!/_{\!\sim}\right)_r$ carries the finest regular topology coarser than the quotient topology on $E\!/_{\!\sim}$, and since
$\left(E\!/_{\!\sim}\right)_s$ was shown to be regular with a topology coarser than the quotient topology, we have continuity 
$\left(E\!/_{\!\sim}\right)_r \to \left(E\!/_{\!\sim}\right)_s \simeq F$.   This means we can apply Theorem \ref{theorem5} again, now with
$G=\left(E\!/_{\!\sim}\right)_r$, and now this implies $\left(E\!/_{\!\sim}\right)_r \simeq \left(E\!/_{\!\sim}\right)_s$, hence the claim $\left(E\!/_{\!\sim}\right)_r\simeq F$.
\hfill $\square$
\end{proof}

\begin{remark}
At first it seems that the correct setting for Corollary \ref{cor1} ought to be semi-regular spaces, not regular spaces.
Unfortunately, the semi-regularization $E_s$ lacks the universal property
of $E_r$, i.e., 
when $f:E \to F$ is continuous and $F=F_s$,
then we do not necessarily get  continuity of $f:E_s \to F$. What is amiss is that
$\tau \subset \tau'$ on $F$ does not imply $\tau_s \subset \tau_s'$. We could also say that the category of semi-regular Hausdorff spaces is not a reflective
subcategory of the category of Hausdorff spaces.
\end{remark}

Yet in the above situation we get $(E\!/_{\!\sim})_s = (E\!/_{\!\sim})_r$, hence $(E\!/_{\!\sim})_s$ is quotient
(extremal epimorphism) in the category of regular spaces, and at the same time has the convenient construction as semi-regularization of $E\!/_{\!\sim}$.

For completely regular spaces we proceed similarly. The complete-regularization of a given space $E$
is the point-set $E$ endowed with the finest completely regular topology  coarser than the given one, denoted $E_{cr}$. Assuming that there exists at least one
Hausdorff completely regular topology coarser than the given one, $E_{cr}$ is Hausdorff.
For $E$ completely regular and  $E\!/_{\!\sim}$ the usual quotient,
we call $(E\!/_{\!\sim})_{cr}$ a completely-regular quotient.  Then a continuous surjection $f:E \to F$ onto a completely regular space $F$ is  completely-regular quotient if
$F \simeq (E\!/_{\!\sim})_{cr}$, where $x\sim y$ iff $f(x)=f(y)$.

\begin{corollary}
\label{cor2}
Let $E, F$ be completely regular and suppose $E$ is a $\tau^*$-Baire space. Let $f:E \to F$ be a continuous nearly open surjection. Then $f$ is 
completely-regular quotient.
Moreover, $F \simeq (E\!/_{\!\sim})_{cr} = (E\!/_{\!\sim})_{s}$. 
\end{corollary}

\begin{remark}
1) The completely-regular quotient may  also be characterized as the initial or limit topology with respect to $C(E\!/_{\!\sim},\mathbb R)$, 
i.e. the coarsest topology on $E\!/_{\!\sim}$ such that all real-valued functions continuous in the quotient topology are continuous.

2)
Yet another way to describe completely-regular quotients is as follows.
Let $\mathscr U, \mathscr V$ be uniformities on $E,F$ inducing the topologies and such that $f$ is uniformly continuous.
(E.g. let $\mathscr U$ be the fine uniformity, then $\mathscr V$ does not matter).
Then the topology on $F$ is the one induced by the quotient uniformity \cite{himmelberg}, and this is $(E\!/_{\!\sim})_s$ under the assumptions of Corollary \ref{cor2}.
\end{remark}

Application of Theorem \ref{theorem4} gives the following analogous:

\begin{theorem}
\label{theorem6}
Let $E$ be a (completely) regular $\tau$-Baire space, $F$ a (completely) regular $\delta$-space. Suppose $f:E \to F$ is a continuous nearly open
surjection. Then $f$ is a (completely) regular-quotient map.
\end{theorem}

\begin{remark}
The reason why we cannot expect $f$ to be open is that quotient maps in the category of Hausdorff spaces need not be open,
and even when they are, $E\!/_{\!\sim}$ need not be semi-regular. One would need openness of $E \to (E\!/_{\!\sim})_s$ to deduce openness of $f$.
\end{remark}

Our third question (c) has now an immediate answer. For $f:E \to F$ continuous nearly open and dense we can expect
$f(E)$ with the topology induced from  $F$  to be a quotient of $E$, and this occurs under the hypotheses of Theorems \ref{theorem5} or \ref{theorem6}.

The above findings motivate the following 

\begin{definition}[Open mapping spaces]
{\rm
Let $\mathscr K$ be a class of Hausdorff spaces. A Hausdorff space $E$ is a $B(\mathscr K)$-space, respectively, a $B_r(\mathscr K)$-space,
if every continuous nearly open surjection, respectively,  bijection, $f$ from $E$ onto any $F \in \mathscr K$ is open.

A Hausdorff space $E$ is a $B_q(\mathscr K)$-space if every
continuous nearly open surjection $f:E\to F$ onto any $F \in \mathscr K$ is a quotient map in the following sense: for
any factorization $f = g\circ h$ with $h:E \to G$ continuous surjective onto a semi-regular Hausdorff space $G$ and
$g:G \to F$ continuous bijective,  it follows that $g$ is a homeomorphism.  When $\mathscr H$ is the class of all Hausdorff spaces, then we say 
$B_q$-space instead of $B_q(\mathscr H)$-space}.
\end{definition}

\begin{remark}
Every $B(\mathscr K)$-space is a $B_q(\mathscr K)$-space, and every semi-regular $B_q(\mathscr K)$-space is a $B_r(\mathscr K)$-space. 
Every $B_r$-space is semi-regular. 
For the latter, observe that the identity $i_E:E \to E_s$ is a continuous nearly open bijection. Since $E_s$ is Hausdorff and $E$ is $B_r$, 
$i_E$ is open, so $E=E_s$. 
\end{remark}

\begin{remark}
Let $\mathscr D$ be the class of $\delta$-spaces. Then Theorem \ref{theorem6} says that every semi-regular $\tau$-Baire space is a 
$B_q(\mathscr D)$-space, hence also a $B_r(\mathscr D)$-space.
\end{remark}

\begin{example}
Let $E = \mathbb R\setminus\mathbb Q \oplus \mathbb R$, $F=\mathbb R$, and let $f:E \to F$ be defined as $f|\mathbb R = i_\mathbb R$,
$f|\mathbb R \setminus \mathbb Q = \iota_{\mathbb R\setminus\mathbb Q}$ the inclusion $\mathbb R\setminus \mathbb Q \to \mathbb R$.
Then $f$ is a continuous nearly open surjection which is not open, as the image of $\mathbb R\setminus \mathbb Q$ is not open in $F$.
But $E$ is completely metrizable, so $f$ is a quotient map, which one can of course see directly as $E\!/_{\!\!\sim} \simeq \mathbb R$ for the
equivalence relation $x\sim y$ iff $f(x)=f(y)$. So $E$ is a $B_q$-space, hence a $B_r$-space, but not a $B$-space. 
\end{example}

\begin{remark}
Recall that a Hausdorff space $E$ is $H$-minimal  \cite[p. 223]{engelking} if there is no strictly
coarser Hausdorff topology on $E$, and it is $H$-closed if it is closed in every Hausdorff space containing $E$ as a subspace. A space is $H$-minimal iff it is $H$-closed and semi-regular; \cite[I.3]{katetov}.  Clearly $H$-minimal spaces are $B_r$-spaces.
\end{remark}

\begin{remark}
In \cite{He} the author constructs a first category $H$-minimal space which according to \cite[p.593]{on_the_theory} 
is not a $B$-space. 
This is in contrast with the following:
\end{remark}

\begin{proposition}
Every
$H$-minimal space is a $B_q$-space. 
\end{proposition}

\begin{proof}
Let $f=g\circ h$ with $h:E \to G$ surjective, $g:G \to F$ bijective, and $G$ semi-regular. 
Since $G$ is the continuous image
of a $H$-closed space $E$, it is $H$-closed by a result
of Kat\v{e}tov; cf. \cite[p. 223]{engelking}, \cite{katetov}. Since $G$ is also semi-regular, it is $H$-minimal.
In consequence $g$ is a homeomorphism.  
\hfill $\square$
\end{proof}

\begin{theorem}
\label{theorem7}
Let $E$ be a semi-regular space containing a dense $B_q$-subspace $D$. Then $E$ is a $B_q$-space. 
\end{theorem}

\begin{proof}
Let $f:E \to F$ be a continuous nearly open surjection, factorized as
$f=g\circ h$ with $h:E\to G$ surjective, $g:G \to F$ bijective, and $G$ semi-regular. We have to show that
$g$ is a homeomorphism.

Since $D$ is dense in $E$, 
Lemma \ref{lemma4} below shows that  the restriction $f|D$ to $D$ remains nearly open as a mapping $D\to F$. But then since
$f(D)$ is dense in $F$,  Lemma \ref{lemma2} shows that
 $f|D$
is still nearly open as a mapping $D\to f(D)$.

Now consider
the restriction  $h|D:D \to h(D)$ of $h$ on $D$, and the restriction $g|h(D): h(D) \to f(D)$ of $g$ on $h(D)$. 
Then $f|D = (g|h(D)) \circ (h|D)$ is a factorization with $g|h(D)$ bijective. Since $h(D)$ is dense in $G$,
it is semi-regular, hence the factorization is amenable to Theorem \ref{theorem5}.
Since
$D$ is by hypothesis a $B_q$-space, $g| h(D)$ is a homeomorphism.

From here, based on
near openness of $g$, denseness of $h(D)$ in $G$, and semi-regularity of $G$, we conclude using
\cite[Lemma p. 589]{on_the_theory}  that $g$ is also a homeomorphism. That completes the argument.
\hfill $\square$
\end{proof}

\begin{lemma}
\label{lemma4}
Let $f:E\to F$ be continuous and surjective. Let $V \subset E$ be open, $D \subset E$ dense.
Then $\xoverline{f(V)} = \xoverline[0.95]{f(V\cap D)}$.
\end{lemma}

\begin{proof}
Let $y\in \xoverline{f(V)}$. Fix a neighborhood $O_F$ of $y$ in $F$. Then
$O_F \cap f(V) \not=\emptyset$, hence
$f^{-1}(O_F)\cap V$ is a nonempty open set in $E$. Since $D$ is dense in $E$, $f^{-1}(O_F)\cap V \cap D \not=\emptyset$, and
that implies $O_F \cap f(V \cap D) \not=\emptyset$. Since $O_F$ was an arbitrary neighborhood of $y$, we have
$y \in \xoverline[0.95]{f(V\cap D)}$.
\hfill $\square$
\end{proof}

\begin{remark}
Theorem \ref{theorem7}
marks a difference between $B$-spaces and $B_q$-spaces. For let $E$ be a non-discrete Lindel\"off $P$-space,
then according to \cite[§7]{on_the_theory} $E$ is a $B$-space and so is the topological sum $E \oplus E$, while
$E \oplus \beta E$ fails to be a $B$-space. Yet, due to the above, $E \oplus \beta E$ is a $B_q$-space.
\end{remark}

\begin{lemma}
\label{concat}
Let $f = g\circ h$ with $g,h$ both continuous and nearly open. Then $f$ is nearly open.
\end{lemma}

\begin{proof}
Let $x\in E$ and $U$ a neighborhood of $x$. Since $h$ is nearly open, we have $h(U) \subset O \subset \xoverline[.9]{h(U)}$ for some open $O$ in $G$.
Since $g$ is also nearly open, we have $f(U)=g(h(U))\subset g(O) \subset \xoverline[.9]{g(O)}^\circ \subset \xoverline[.9]{g(O)} \subset \xoverline[1]{g(\xoverline[.9]{h(U)})}
\subset \xoverline[1]{g(h(U))} = \xoverline[.9]{f(U)}$ using continuity of $g$.
\hfill $\square$
\end{proof}

For this result see also \cite[Lemma 8]{tkachenko}.

\begin{proposition}
Let $E$ be a $B_q$-space and $h$ a continuous nearly open surjection onto a semi-regular space $G$. Then $G$ is a $B_r$-space.
\end{proposition}

\begin{proof}
Let $g:G \to F$ be a continuous nearly open bijection onto the space $F$. Then $f=g\circ h:E \to F$ is a continuous surjection, which by Lemma 
\ref{concat} is nearly open. As $E$ is a $B_q$-space, $g$ is a homeomorphism. 
That proves the claim.
\hfill $\square$
\end{proof}

\section{Determined game}
\label{sect_determined}
An infinite two-person game with two possible outcomes is called {\it determined} if either player $\alpha$
or player $\beta$ has a winning strategy. The classical BM-game is not determined, as there exist Baire spaces which are not
weakly $\alpha$-favorable. In contrast, we shall see that for certain classes of spaces $E$ the ${\rm BM}'$-game is determined.

\begin{proposition}
\label{prop3}
Let $E$ be a metrizable locally convex vector space which is $\tau$-Baire. Then $E$ is complete.
\end{proposition}

\begin{proof}
Since $E$ is a $B_r$-space by Theorem \ref{theorem5}, it is also a $B_r$-complete
separated locally convex vector spaces in the sense of Pt\'ak. Hence, as the nomenclature suggests,  $E$ is complete, 
as follows from
\cite[§34, 2.(1)]{koethe}.
\hfill $\square$ 
\end{proof}

This means the ${\rm BM}'$-game is determined in the class of metrizable locally convex vector spaces $E$. 
Because if $\beta'$ has no winning strategy, then $E$ is $\tau$-Baire, hence is completely metrizable by Proposition \ref{prop3}, so that
by Proposition  \ref{prop1} player $\alpha'$ has a winning strategy.
The same is true when we work in the category of separated topological vector spaces, where
$B$- and $B_r$-completeness can be defined accordingly \cite{adasch}.

\begin{proposition}
\label{prop4}
Let $G$ be a metrizable topological group which is $\tau$-Baire. Then
$G$ is completely metrizable.
\end{proposition}

\begin{proof}
Let $\widetilde{G}$ be the completion  of $G$ in its two-sided uniformity. Suppose $\widetilde{x}\in \widetilde{G}\setminus G$.
We form the space $F = G \cup \widetilde{x}G$ equipped with the topology induced from $\widetilde{G}$. Let $G \oplus G$ be the 
topological sum in the sense of \cite[p. 74]{engelking} and 
define
the function $f: G \oplus G \to F$ as follows. Let the point-set
of $G \oplus G$ be $G\times \{1\} \cup G \times \{2\}$, and let $f(x,1) = x$, $f(x,2)=\widetilde{x}x$. Then $f$ is continuous and bijective,
the latter since $G \cap \widetilde{x}G =\emptyset$. Since $G$ and $\widetilde{x}G$ are both dense in $\widetilde{G}$,
the mapping $f$ is nearly open by Lemma \ref{lemma4}. Since $G$ is $\tau$-Baire and metrizable, it is a $B_r$-space. By \cite[Thm. 1]{sums_and_products} the topological sum
$G \oplus G$  is also $B_r$, hence $f$ is open, so that both $G$ and $\widetilde{x}G$ are open in $\widetilde{G}$.
That, however, is impossible. Hence $G = \widetilde{G}$.
\hfill $\square$
\end{proof}

\begin{remark}
The reasoning here differs from the one in Proposition \ref{prop3}, because $B_r$-groups  need not be complete. 
For instance $Q=\{e^{2\pi iq}: q\in  \mathbb Q\}$ is a $H$-minimal group, \cite{Gr}, hence a $B_r$-group. Clearly $Q$ is not a $B_r$-space in the topological sense.
\end{remark}

\begin{remark}
Proposition \ref{prop4} shows that in a metrizable topological group, one of the players $\alpha',\beta'$ has a winning strategy, so again in that category the
${\rm BM}'$-game is determined.
\end{remark}

\begin{remark}
\label{harder}
In Remark \ref{easy} we remarked that Baire spaces need not be $\tau$-Baire. Now we remark that even spaces $E$ where $E\times E$
is Baire need not be $\tau$-Baire. Let $E$ be a separable metrizable locally convex vector space which is Baire but not
complete. Then $E \times E$ is Baire by a result of Kuratowski \cite[p. 201]{engelking}, but by Proposition \ref{prop3} above $\beta'$
has a winning strategy, so $E$ is not $\tau$-Baire. This also shows that  the converse of Proposition \ref{beta1} is incorrect.
\end{remark}

\begin{remark}
Kuratowski's result gives even $E^n$  Baire for every $n$.
\end{remark}

Let $G$ be a topological group acting continuously on a Hausdorff space $X$. The $G$-flow on $X$ is said to be minimal if every orbit $Gx$
is dense in $X$. Then, expanding on the technique in Proposition \ref{prop4}, we have the following

\begin{proposition}
Suppose the $G$-flow is minimal. 
Let $E = Gx$ be an orbit which is a $B_r$-space in the topology induced by $X$. Then $Gx$ is not homeomorphic to any other orbit $Gx'$.
\end{proposition}

\begin {proof}
Let $E = \{gx: g\in G\}$ for some $x\in X$, and suppose $E \simeq E' = \{gx': g\in G\}$ for another $x'\in X$. Let $h:E \to E'$ be a homeomorphism.
We let $E \oplus E$ be the topological sum
of two copies of $E$ represented on the set $E \times \{1\} \cup E \times \{2\}$. Then we define
$f: E \oplus E \to X$ as $f(x,1) = x$, $f(x,2) = h(x)$. Then $f$ is a continuous bijection onto the dense subset $E \cup E'$ of $X$, where injectivity is due to
$E \cap E'=\emptyset$. Since
$E \xhookrightarrow{} X$ and $E' \xhookrightarrow{} X$ are dense embeddings and $h$ is a homeomorphism, $f$ is nearly open. Since $E$ is a $B_r$-space, so is $E \oplus E$. Hence
$f$ is a homeomorphism onto $E \cup E'$. But  that means $E,E'$ are both open in $E \cup E'$, contradicting the fact that they are both dense.
\hfill $\square$ 
\end{proof}

\begin{remark}
It is known \cite{yaacov} that if $G$ is a Polish group, $X$ a compact Hausdorff space, and if the $G$-flow on $X$ is minimal and metrizable, then
there exists  an orbit $Gx_0$ which is residual in $X$. That means $Gx_0$ contains a dense $G_\delta$ subset of $X$, hence has a dense 
\v{C}ech-complete subspace, hence is a $B_r$-space \cite{BP,on_the_theory}. Then $Gx_0$ cannot be homeomorphic to any other orbit.  

However, even when no residual orbit exists, 
as soon as we have a $\tau^*$-Baire orbit, it is not homeomorphic to any other orbit. 
\end{remark}

\section{Barely Baire spaces}
\label{sect_barely}
We recall a construction from \cite{fleissner}. 
A subset $S \subset \omega_1$ is cofinal if $|S|=\omega_1$,
and it is closed when closed in the order topology on $\omega_1$. A subset $S \subset \omega_1$
is stationary if it intersects every closed cofinal subset of $\omega_1$. Let $\omega_1$ be endowed with the discrete topology and
give $\omega_1^\omega$ the product topology. Then $\omega_1^\omega$ is metrizable. For $S \subset \omega_1$ cofinal put
$S^*=\{f\in \omega_1^\omega: f^*:=\sup_n f(n) \in S\}$, endowed with the topology induced by the product topology. Then by
\cite[Example 1, p. 234]{fleissner} $S^*$ is a Baire space iff $S$ is stationary. Moreover, if $S,T$ are two stationary subsets
of $\omega_1$ such that $S \cap T$ is not stationary, then $S^*,T^*$ are Baire, but $S^* \times T^*$ is not, hence  $S^*,T^*$ are 
what the authors of \cite{fleissner} call  {\it barely Baire}.
In addition, if we arrange $S \cap T=\emptyset$, then the topological sum $S^*\oplus T^*$
is no longer a $B_r$-space, nor is the product $S^* \times T^*$, see \cite[p. 183]{fleissner}.
In \cite{fleissner} it is also proved that $S^*$ contains a dense completely metrizable subspace iff $S$ is
closed cofinal, which here due to metrizability is equivalent to weak $\alpha$-favorability of $S^*$, or to $p$- and $c$-completeness. 

\begin{proposition}
\label{prop9}
Let $S$ be stationary. Then $S^*$ is a $\tau$-Baire space, i.e., player $\beta'$ has no winning strategy.
\end{proposition} 

\begin{proof}
Recall that sets of the form $B(\gamma) = \{f\in S^*: f(i) = \gamma_i, i=1,\dots,k\}$ for finite sequences
$\gamma=(\gamma_1,\dots,\gamma_k)$ of ordinals $\gamma_i<\omega_1$ form a basis for $S^*$. 
We consider a strategy $\beta'$ playing with basic sets.
Suppose $V=B(\gamma)$ and $V'=B(\gamma')$ are basic open sets with $V' \subset V$,  then $\gamma'$ is a prefix of $\gamma$. 
Therefore play sequences
$V_1,V_1',W_1,W_1',\dots,V_r,V_r',\dots$ 
give rise to finite sequences of ordinals  $\gamma_1 \subset \gamma_1' \subset \gamma_2 \subset \gamma_2' \subset \dots$ and 
$\delta_1 \subset \delta_1' \subset \delta_2\subset \delta_2'\subset \dots$. 
We may therefore consider $\beta'$ as a function on finite sequences of ordinals
$\beta'(\gamma_1,\gamma_1',\delta_1,\delta_1',\dots,\delta_r,\delta_r') = \gamma_{r+1}$,
respectively $\beta'(\gamma_1,\gamma_1',\delta_1,\delta_1',\dots,\gamma_r,\gamma_r') = \delta_{r}$.

We say that an ordinal $\eta < \omega_1$ is a fixed point of $\beta'$
if  $\max(\gamma_r) < \eta$, $\max(\delta_r) < \eta$ together with
$\beta'(\gamma_1,\gamma_1',\delta_1,\delta_1',\dots,\delta_r,\delta_r') = \gamma_{r+1}$ imply
$\max(\gamma_{r+1}) < \eta$, and at the same time $\max(\delta_{r-1}) < \eta$, $\max(\gamma_r) < \eta$ together with
$\beta'(\gamma_1,\gamma_1',\delta_1,\delta_1',\dots,\gamma_r,\gamma_r') = \delta_{r}$, imply
$\max(\delta_r) < \eta$. 

We let
$F$ be the set of fixed points of $\beta'$. Clearly $F$ is cofinal and closed. 
Now using stationarity of $S$, let $\eta\in S \cap F$.
Fix a sequence $\eta_1 < \eta_2 < \dots$ converging to $\eta$. Now define
$\alpha'(\gamma_1,\gamma_1',\delta_1,\delta_1',\dots,\gamma_r) = \gamma_r^{\frown} \eta_r =:\gamma_r'$
and $\alpha'(\gamma_1,\gamma_1',\delta_1,\delta_1',\dots,\delta_r) = \delta_r^{\frown} \eta_r =:\delta_r'$. Then the play of $\alpha'$
against $\beta'$
generates two sequences $\delta_1 \subset \delta_1' \subset \delta_2\subset \delta_2' \subset \dots$
and  $\gamma_1 \subset \gamma_1' \subset \gamma_2\subset \gamma_2' \subset \dots$ with $\sup_i \gamma_i = \sup_i \delta_i = \sup_i\eta_i = \eta$.
Let $f = \bigcup_{i=1}^\infty \gamma_i$ and $g=\bigcup_{i=1}^\infty \delta_i$, then
$f^*=g^*=\eta\in S$, hence $f,g\in S^*$. We have $f\in \bigcap_{k=1}^\infty B(\gamma_k)$ and $g\in \bigcap_{k=1}^\infty B(\delta_k)$,
hence $\alpha'$ wins against $\beta'$.
\hfill $\square$
\end{proof}

\begin{corollary}
If $S$ is stationary, then $S^* \times S^*$ is Baire, $S^*$ is a $B_q$-space, and hence a $B_r$-space.
\end{corollary}

\begin{proof}
The first statement is from \cite{fleissner}. The second statement follows with Proposition \ref{prop9} 
in tandem with Theorem \ref{theorem6}.
\hfill $\square$
\end{proof}

\begin{corollary}
Suppose $S$ is stationary, but does not contain any closed cofinal set. Then neither
$\alpha'$ nor $\beta'$ have winning strategies, so in this class of spaces $S^*$ the ${\rm BM}'$-game is not determined.
\end{corollary}

\begin{proof}
By Proposition \ref{prop9} player $\beta'$ has no winning strategy. 

On the other hand, suppose $\alpha'$ has a winning strategy on $S^*$. Then by Proposition \ref{prop1}, so has $\alpha$.
As in the proof of Proposition \ref{prop9} we can consider $\alpha$ 
a function on finite sequences of ordinals: $\alpha(\gamma_1,\delta_1,\gamma_2,\delta_2,\dots,\gamma_r) = \delta_r$, where
$\gamma_1 \subset \delta_1 \subset \gamma_2 \subset \delta_2 \subset\dots$. 
Let $F$ be the set of fixed-points of $\alpha$, then $F$ is closed cofinal in $\omega_1$.  We show that $F \subset S$.
Indeed, let $\eta \in F$ and fix a sequence $\eta_1 < \eta_2 < \dots$ converging to $\eta$. Define a strategy $\beta$ for the second player as follows:
$\beta(\gamma_1,\delta_1,\dots,\gamma_r,\delta_r) = \delta_r^{\frown}\eta_r =: \gamma_{r+1}.$ Let $\gamma_1,\delta_1,\dots$ be the play of $\alpha$ against $\beta$. 
By the definition of $F$ and the $\eta_i$ we have $\sup_i \gamma_i = \sup_i \delta_i = \eta$. But $\alpha$ wins against $\beta$, so 
the function $f=\bigcup_i \gamma_i = \bigcup_i \delta_i$ belongs to $S^*$, which implies $\eta = f^*\in S$. That proves the second claim by contraposition.
\hfill $\square$
\end{proof}

\begin{remark}
Amusingly, if we consider only cofinal sets $S \subset \omega_1$ which are Borel in the order topology of $\omega_1$, then
either $S$ or $\omega_1\setminus S$ contains a closed cofinal set. Hence in this class of spaces $S^*$, the ${\rm BM}'$-game is again determined. 
\end{remark}

\section{Barely $\tau$-Baire spaces}
\label{sect_barely_tau}
We recall a second construction from \cite{fleissner}. Let $c=2^\omega$, $c^+$ the successor
cardinal of $c$, $C_\omega c^+$ the set of ordinals $\alpha < c^+$ with cofinality cf$(\alpha) \leq \omega$. 
We select a family $\{A_y: y\in 3^\omega\}$ of mutually disjoint stationary subsets of $C_\omega c^+$. Then
for $x\in 3^\omega$ define
$$
B_x = \bigcup \{A_y: y \in 3^\omega, y(n) \not= x(n) \mbox{ for all $n\in \omega$}\}.
$$
Note that this is in particular possible because $C_\omega c^+$ is itself a stationary
subset of $c^+$. Now we endow $3$ and $c^+$ with the discrete topology  and  define
$$
E = \{(x,f) \in 3^\omega \times (c^+)^\omega: f^* = \sup_n f(n) \in B_x\},
$$
endowed with the  topology induced by the product topology. Since $3^\omega \times (c^+)^\omega$ is metrizable, so is $E$.
Basic sets of $E$ are
indexed by finite sequences $\sigma = ((\alpha_0,\beta_0),\dots,(\alpha_r,\beta_r))\in (3\times (c^+))^{r+1}$,
that is,
$$
B(\sigma)=\{(x,f) \in E: x(i)=\alpha_i, f(i) = \beta_i, i=0,\dots,r\}.
$$

Following the line of \cite[Example 3, p. 675]{baire_category} we have

\begin{proposition}
The space $E$ is $\tau$-Baire.
\end{proposition}

\begin{proof}
Let $\beta'$ be a strategy  in the  ${\rm BM}'$-game. We may assume that it plays with basic open sets 
$B(\sigma)$, $B(\tau)$. Again if $B(\sigma) \subset B(\sigma')$, then $\sigma' \subset \sigma$ is sequence extension.
A play gives therefore rise to two sequences $\sigma_1 \subset \sigma_1' \subset \sigma_2 \subset \sigma_2' \subset \dots$
and $\tau_1 \subset \tau_1' \subset \tau_2 \subset \tau_2' \subset \dots$ intertwined as in the figure.  
Hence a strategy $\beta'$ defines a mapping, also denoted $\beta'$, with
$$
\beta'(\sigma_1,\sigma_1',\tau_1,\tau_1',\dots,\sigma_r,\sigma_r') = \tau_r,
\quad \beta'(\sigma_1,\sigma_1',\dots,\tau_r,\tau_r') = \sigma_{r+1}.
$$
We can see this as two mappings $\Theta_1,\Theta_2$ on finite sequences $\sigma = (\sigma_1,\sigma_1',\sigma_2,\sigma_2',\dots)$ and 
$\tau=(\tau_1,\tau_1',\tau_2,\tau_2',\dots)$ such that
$\Theta_1(\sigma,\tau) \supset \tau$ and $\Theta_2(\sigma,\tau) \supset \sigma$. 
Now let $W$ be the set of ordinals $2 < \alpha < c^+$ with the following property:

{\em
There exist sequences $(\sigma_n)$ and $(\tau_n)$ such that $\Theta_1(\sigma_0,\tau_0) = \sigma_1$, $\sigma_1 \subset \sigma_1'$,
$|\sigma_1'|=|\sigma_1|+1$, $\Theta_2(\sigma_1,\tau_0) = \tau_1$, $\tau_1 \subset \tau_1'$, $|\tau_1'| = |\tau_1|+1$, etc., such that
$\bigcup_n\sigma_n = (x,g)$, $\bigcup_n\tau_n =(y,h)$ with $g^*=h^*=\alpha$.
}

The set $W$ is stationary in $c^+$. Indeed, if $C$ is closed cofinal in $c^+$, define the extensions $\sigma_i'$ of the $\sigma_i$
and $\tau_i'$ of the $\tau_i$ such that $\max \sigma_i' \leq \max \tau_{i+1}' \leq \max \sigma_{i+1}'$ and such that 
$\sup_i \sigma_i' = \sup_i \tau_i' \in C$. Then
$\bigcup_n \sigma_n =(x,g)$ and $\bigcup_n \tau_n = (y,h)$ have $g^*=h^* \in C \cap W$. This proves stationarity of $W$.

Now fix $z\in 3^\omega$ and let $W_z$ be the set of those $\alpha\in W$ where
$(x,g)$, $(y,h)$ exist as above, with $g^*=h^*=\alpha$, but $z(n)\not\in \{x(n),y(n)\}$ for all $n$. Then
$$
W = \bigcup\{W_z: z\in 3^\omega\}
$$
and since $W$ cannot be the union of $c < c^+$ stationary sets, one of the $W_z$ is stationary. Now we use the following auxiliary result
\cite[Lemma 1]{fleissner}:

{\em
If $K \subset (c^+)^\omega$ is closed and $W=\{f^*: f \in K\}$ is stationary, then there exists
a closed cofinal subset $C$ of $c^+$ such that
$C \cap C_\omega c^+\subset W$.
}

We
apply this to the stationary $W_z$ found above in order to find $C \cap C_\omega c^+ \subset W_z$. For that let
$K$ be the set of all $((x,g),(y,h))$ with $\bigcup_n \sigma_n = (x,g)$ and $\bigcup_n \tau_n=(y,h)$, $g^*=h^*=\alpha$,
for sequences as in the definition of $W$, having
$z(n)\not\in \{x(n),y(n)\}$ for every $n$. Then $K$  is closed in $3^\omega \times (c^+)^\omega \times 3^\omega \times (c^+)^\omega$, and 
we have $W_z=\{\psi^*: \psi\in K\}$.

As the sets $A_y$ arising in the construction of $E$ are stationary, we have $W_z \cap A_z \cap C_\omega c^+\not=\emptyset$. Choose
$\gamma$ herein and let $(\sigma_n)$, $(\tau_n)$ be sequences as above
giving rise to $\bigcup_n \sigma_n=(x,g)$, $\bigcup_n\tau_n = (y,h)$ with $g^*=h^*=\gamma$ and $z(n)\not\in \{x(n),y(n)\}$.
Since $\gamma \in A_z$, we have $(x,g),(y,h)\in E$. But $(x,g)$ is in the intersection of the $B(\sigma)$,
and $(y,h)$ is in the intersection of the $B(\tau)$, hence both intersections are non-empty, and strategy $\alpha'$ is winning against
strategy $\beta'$. That proves the claim.
\hfill $\square$ 
\end{proof}

\begin{remark}
The interest in this space is that $E\times E$ is Baire, (cf. Theorem \ref{theorem3}), 
but $E\times E \times E$ is no longer Baire \cite{fleissner}. This gives a metrizable $\tau$-Baire space $E$ whose square $E\times E$ is
no longer $\tau$-Baire (even though $E\times E$ is Baire).

Secondly this gives rise to a metrizable $B_r$-space whose square is no longer $B_r$.
Namely, either  (a) $E\times E$ is not a $B_r$-space. Then $E$, which is $B_r$ by Theorem \ref{theorem5a},  is the space we are looking for. 
Or (b) $E\times E$ is again $B_r$. Then $(E\times E)\times(E\times E)$ is no longer a $B_r$-space,
because it is not Baire, but due to \cite{baire_category} ought to be Baire if it were $B_r$.
So here we have the $B_r$-space $E\times E$ whose square is no longer $B_r$. 
\end{remark}

\begin{remark}
It would be interesting to know which of the two cases (a), (b) above is true. 
If (b) holds with $E \times E$  still $B_r$, then we have
a metrizable $B_r$-space $E \times E$, which is Baire, but not $\tau$-Baire.  Currently no such space is known.
\end{remark}

\section{Michael game}
\label{sect_michael}
The following variation of the Banach-Mazur game was introduced by Michael \cite{michael_game}. Players
$\beta$ and $\alpha$ choose successively non-empty sets $B_1 \supset A_1 \supset B_2 \supset A_2 \supset \dots$ such that
$A_i$ is open in $B_i$, that is, $A_i = B_i \cap U_i$ for some open set $U_i$ in $E$. Player $\alpha$ playing with the sets $A_i$
wins the game if $\bigcap_{n=1}^\infty \closure[2]{A}_n \not=\emptyset$. Player $\alpha$ wins strongly
if every filter $\mathscr F$ with $A_i\in \mathscr F$ for every $i$ has a cluster point, i.e.,
$\bigcap \{\closure[2]{F}: F \in \mathscr F\} \not=\emptyset$.
 In \cite[Thm. 7.3]{michael_game} the author proves that $\alpha$ has a winning strategy
 if and only if $E$ has a complete, exhaustive sieve. 
 
 We say that $E$ is a $m$-Baire space if player $\beta$ does not have a winning strategy in the Michael game.
 It is clear that every $m$-Baire space is Baire, because if player $\beta$ plays with open sets, then
 $\alpha$ automatically responds with open sets, and  the play coincides with the Banach-Mazur game. 
 
 \begin{proposition}
 Let $E$ be a regular $m$-Baire space.
 Then every $G_\delta$ subset $G$ of $E$ is a $m$-Baire space.
 \end{proposition}
 
 \begin{proof}
 Let $G = \bigcap_{n=1}^\infty G_n$ with open sets $G_1 \supset G_2 \supset \dots$, and let 
$\beta_G$ be a strategy for player $\beta$ in the Michael game on $G$. We define a strategy $\beta$ on the whole space.

Suppose $\beta_G(\emptyset) = B_1$. By regularity we may choose a set $V_1$ open in $E$ such that $\overline{V}_1 \subset G_1$ and $B_1' = B_1 \cap V_1 \not=\emptyset$.
Then define $\beta(\emptyset) = B_1'$. Now let $A_1 \subset B_1'$ be non-empty and relatively open in $B_1'$, i.e., $A_1 = B_1' \cap U_1$ for some open $U_1$
in $E$. We have to define $\beta(B_1',A_1)$. Note that $A_1 \subset B_1' \subset B_1$ and
$A_1 = B_1' \cap U_1 = B_1 \cap V_1 \cap U_1$, and since $A_1 \subset G$, $A_1 = B_1 \cap (V_1 \cap U_1 \cap G)$, so 
$A_1$ is relatively open in $B_1$ with regard to the space $G$. Hence $B_2 = \beta_G(B_1,A_1)$ is defined and is a non-empty subset of $A_1$.
We choose a set $V_2$ open in $E$ such that $\overline{V}_2 \subset G_2$ and $B_2' = B_2 \cap V_2 \not=\emptyset$.
Define $\beta(B_1',A_1)=B_2'$. Etc.

By assumption $E$ is $m$-Baire, so there exists a strategy $\alpha$ winning against $\beta$. Let us define a
strategy $\alpha_G$ on $G$ which wins against $\beta_G$. We have to define
$\alpha_G(B_1)$, where $B_1 = \beta_G(\emptyset)$. We recall the construction of $B_1' \subset B_1$, $B_1' = \beta(\emptyset)$ and
define $\alpha_G(B_1) = \alpha(B_1') = A_1$. Next we have to define
$\alpha_G(B_1,A_1,B_2)$. We recall the construction  of $B_2' \subset B_2$ and define
$\alpha_G(B_1,A_1,B_2) = \alpha(B_1',A_1,B_2') = A_2$. Etc. It is clear from the construction that
the play of $\alpha_G$ against $\beta_G$ and the play of $\alpha$ against $\beta$ are interlaced  as follows:
$B_1 \supset B_1' \supset A_1 \supset B_2' \supset B_2 \supset A_2 \supset \dots$. Since $\alpha$ wins, we have
$\bigcap_{n=1}^\infty \closure[2]{A}_n \not=\emptyset$, where the closure is with respect to $E$. Choose
$x$ within. Since $\closure[2]{A}_i \subset \overline{B_i \cap V_i}\subset \overline{V }_i \subset G_i$, we have $x\in G$.
But $A_i \subset G$, hence $cl_G(A_i) = \overline{A}_i \cap G$, hence $x\in cl_G(A_i)$. That shows
$\bigcap_{n=1}^\infty cl_G(A_n) \not=\emptyset$, hence $\alpha_G$ wins against $\beta_G$.
\hfill $\square$
 \end{proof}
 
 Since in a metrizable space closed sets are $G_\delta$-sets, we have the following
 
 \begin{corollary}
 Every metrizable m-Baire space $E$ is hereditary Baire. In particular, $E \times F$ is Baire for every Baire space $F$.
 \end{corollary}

\begin{proof}
This follows from a result of Moors \cite{moors_hereditary_baire}.
\hfill $\square$
\end{proof} 

\begin{remark}
As a consequence we see that the $\tau$-Baire spaces $S^*$ of Section \ref{sect_barely} and \cite{fleissner} are not $m$-Baire, because if $S,T$ are disjoint
stationary subsets of $\omega_1$, then $S^*,T^*$ are metrizable Baire, but $S^*\times T^*$ is not Baire. Hence $S^*$ is not hereditary Baire by  
\cite{moors_hereditary_baire}, and so
cannot be $m$-Baire. 
 \end{remark}
 
In yet another well-known modification of the Banach-Mazur game
player $\beta$ chooses open sets $V_k$ and points $x_k \in V_k$, while player $\alpha$ has to respond with open sets
$U_k$ satisfying $x_k \in U_k \subseteq V_k$. Player $\beta$ wins when $\bigcap_{k=1}^\infty U_k = \emptyset$, otherwise $\alpha$ wins. 
The game is played with perfect information. The space is said to be strongly $\alpha$-favorable if player $\alpha$ has a winning strategy
\cite{debs}. 
We are interested in the undetermined case, i.e. when neither $\alpha$ nor $\beta$ have winning strategies. This game was proposed by Choquet
\cite{choquet} and is thoroughly discussed in \cite{debs}, and we
write $\alpha^c,\beta^c$ for the corresponding 
strategies, all assumed to have complete memory. In the $\beta^c$-defavorable case we call the space $c$-Baire.

\begin{proposition}
\label{prop13}
Every metrizable $c$-Baire space is $m$-Baire. 
\end{proposition}

\begin{proof}
Let $E$ be metrizable and
consider a strategy $\beta^m$ in the Michael game. We define an associated strategy $\beta^c$.
Suppose $\beta^m(\emptyset) = B_1$. We pick $b_1\in B_1$ and choose an open set $V_1$ with $d(V_1)\leq 1/1$ and $b_1 \in V_1$.
Put $\beta^c(\emptyset)=(b_1,V_1)$.
Suppose now $U_1$ satisfies $b_1 \in U_1 \subset V_1$, so that it could be a move of $\alpha$ in response to the move
$\beta^c(\emptyset)$. In that case we put $A_1 = B_1 \cap U_1$ and interpret $A_1$ as a move of $\alpha$ in response to
$\beta^m(\emptyset)$.

Suppose now that $B_2=\beta^m(B_1,A_1)$. If $A_1 = B_1 \cap U_1$ with $U_1$ constructed as above, then
$V_1$ is available due to perfect memory, and
we want to define $\beta^c((b_1,V_1),U_1) = (b_2,V_2)$. We pick $b_2 \in B_2$ and $V_2$ such that $b_2 \in V_2 \subset U_1$
satisfying $d(V_2)\leq 1/2$. Then $(b_2,V_2)$ is our $\beta^c$ move.

Continuing in this way gives a strategy $\beta^c$ in the Choquet game. Let $\alpha^c$ be a strategy
which wins against $\beta^c$. We use it to define a strategy $\alpha^m$ which wins against $\beta^m$. 
We have to define $\alpha^m(B_1,A_1,\dots,B_i)$, where $A_j = B_j \cap U_j$ for $1 \leq j \leq i-1$. In those cases where
$\beta^c$ has been derived from $\beta^m$ as above, we have access to the corresponding $(b_i,V_i)$, so we get
$\alpha^c((b_1,V_1),U_1,(b_2,V_2),U_2,\dots,(b_i,V_i)) = U_i$. We then define $\alpha^m(B_1,A_1,\dots,B_i) = A_i := B_i \cap U_i$.
In all other cases we define $\alpha^m$ at leisure.

Suppose now $B_1,A_1,B_2,A_2,\dots$ is the play of $\alpha^m$ against $\beta^m$. Then
we get the play $(b_1,V_1),U_1,(b_2,V_2),U_2,\dots$ of $\alpha^c$ against $\beta^c$, and here $\alpha^c$ is winning, so
we have $\bigcap_{i=1}^\infty U_i\not=\emptyset$. Pick $\bar{x}$ herein, then due to $d(V_i)\leq 1/i$ we have $b_i \to \bar{x}$, hence
$\bar{x} \in \bigcap_{i=1}^\infty \closure[2]{B}_i=\bigcap_{i=1}^\infty \closure[2]{A}_i$, hence $\alpha^m$ wins against $\beta^m$.
\hfill $\square$
\end{proof}

\begin{remark}
\label{remark16}
As a consequence, $\tau$-Baire spaces $S^*$ are not $c$-Baire either, i.e., for $S\subset \omega_1$ stationary but not containing a closed cofinal
subset, $S^*$ is $\beta^c$-favorable, even though a $\tau$-Baire space.
\end{remark}
 
 \begin{remark}
 A variant of the BM-game which bears some resemblance with our ${\rm BM}'$-game is the Reznichenko-game of \cite{moors_semi},
 but neither game seems stronger than the other. 
 \end{remark}

\section{Closed Graph Theorem}
\label{closed_graph}
Moors \cite[Thm. 2]{moors1} proves that every nearly continuous closed graph mapping $f:E \to F$ from a
Baire space $E$ to a partition complete space $F$ is continuous. Partition complete spaces, also known as cover complete spaces
\cite{michael_game}, are those where player $\alpha$ has a strong winning strategy in the Michael game.

The key observation here  is that
$c$-completeness of $F$ is by a little margin too weak to prove 
the closed graph theorem, 
which is why in
\cite{graph_theorem} a notion called strict $c$-completeness, equivalent to 
$\alpha$-favorability in the strong Choquet game, was used.  But that notion  is now by a little margin too strong.  
The point made by
\cite{moors1} is that partition completeness, settled in between these twain, is just about right.  
This can also be seen in the light of Proposition \ref{prop13} and Remark \ref{remark16}.

There are two ways to expand from here. We may introduce a tandem Michael game on $F$ and weaken
$\alpha$-favorability to $\beta'$-defavorability, while strengthening Baire category of $E$ to $\alpha$-favorability.
A second option is to keep the weaker hypotheses: Baire category  of $E$ and $c$-completeness of $F$, and require 
instead a little more on $f$. 
We shall 
follow this second line.

\begin{lemma}
\label{help}
Let $f:E \to F$ be nearly open and nearly continuous  with closed graph. Let $E$ be Baire and
$F$  $c$-complete. Suppose $G$ is a dense $G_\delta$ in $F$, $V$ open in $F$, $H$ a dense $G_\delta$ in $E$.
Then $\xoverline[.92]{f^{-1}(V)} = \xoverline[1.0]{f^{-1}(V\cap G) \cap H}$.
\end{lemma}

\begin{proof}
It suffices to prove
$f^{-1}(V) \subset \xoverline[1.0]{f^{-1}(V\cap G) \cap H}$. Let $z\in f^{-1}(V)$, and take an open neighborhood
$W$ of $z$. We have to show $W \cap f^{-1}(V\cap G) \cap H\not=\emptyset$.

Let $(T,\phi)$ be a web on $F$ satisfying $(w_1)$, $(w_2)$ and $(c)$ from Section \ref{prepare}. 
Write $H = \bigcap_{n=1}^\infty H_n$ and $G=\bigcap_{n=1}^\infty G_n$
with $G_n,H_n$ dense open and decreasing. We shall define a strategy $\beta$
for the BM-game on $E$.

We have to define $\beta(\emptyset)$.
The set $X_1= \bigcup \{\phi(t): t\in T, \phi(t) \subset V \cap G_1\}$ is dense in $V$, while $W \cap H_1$ is dense in $W$.
Now $\xoverline[.9]{f(\xoverline[0.95]{W\cap H_1})} \cap V$ is a neighborhood of $f(z)$, hence it intersects $X_1$, as $X_1$ is dense in $V$.
By the definition of $X_1$ there exists $x_1$ in this intersection and $t_1\in T$ with $x_1 \in \phi(t_1) \subset V \cap G_1$. 
Now $\phi(t_1)\cap f(\xoverline[.95]{W \cap H_1}) \not=\emptyset$. Choose $y_1 \in \xoverline[1.0]{W\cap H_1}$ with $f(y_1)\in \phi(t_1)$. 
Since $\xoverline[.9]{f^{-1}(\phi(t_1))}$ is a neighborhood of $y_1$,  it intersects $W\cap H_1$. Let $z_1$ be in this intersection. We pick an open set 
$W_1$ with $z_1 \in W_1 \subset W \cap H_1 \cap\xoverline[.9]{f^{-1}(\phi(t_1))}$. There exists $w_1 \in W_1 \cap f^{-1}(\phi(t_1))$.
Let $\beta(\emptyset) = W_1$.

Let $W_1' \subset W_1$ be nonempty open, then we have to define $\beta(W_1,W_1')$. Now 
$D_2 = W_1' \cap H_2$ is dense in $W_1'$, while $X_2 = \bigcup\{\phi(t): t_1 <_Tt, \phi(t) \subset G_2\}$ is dense in $\phi(t_1)$. 
Hence $\xoverline[.9]{f(\xoverline[.9]{D}_2)} \cap X_2$ is nonempty. By the definition of $X_2$ there exists $x_2$ in this intersection
and $t_2$ with $t_1 <_Tt_2$ and $x_2 \in \phi(t_2) \subset G_2$.
That implies $\phi(t_2) \cap f(\xoverline[.9]{D}_2) \not=\emptyset$. Choose $y_2 \in \xoverline[.9]{D}_2$ with $f(y_2)\in \phi(t_2)$.
Then $\xoverline[.9]{f^{-1}(\phi(t_2))}$ is a neighborhood of $y_2$, so cuts $D_2$. Choose $z_2\in \xoverline[.9]{f^{-1}(\phi(t_2))}\cap D_2$
and an open $W_2$ with $z_2 \in W_2 \subset W_1' \cap H_2 \cap  \xoverline[.9]{f^{-1}(\phi(t_2))}$. There exists $w_2 \in W_2 \cap f^{-1}(\phi(t_2))$.
We let $\beta(W_1,W_1')=W_2$.

Continuing in this way defines a strategy $\beta$ in the ${\rm BM}$-game. Let $\alpha$ be a strategy winning against $\beta$.
Let $W_1\supset W_1' \supset W_2 \supset W_2'\supset \dots$ be their play. Then by construction  we have
\begin{itemize}
\item[i.] $D_i = W_{i-1}' \cap H_i$ dense in $W_{i-1}'$, 
\item[ii.] $X_i = \bigcup\{\phi(t): t_{i-1} <_Tt, \phi(t) \subset G_i\}$ dense in $\phi(t_{i-1})$, $x_i \in \phi(t_i) \subset G_i$, $t_{i-1} <_Tt_i$.
\item[iii.] $y_i \in \xoverline[1.0]{W_{i-1}' \cap H_i}$, $z_i \in W_i \subset W_{i-1}' \cap H_i \cap \xoverline[.9]{f^{-1}(\phi(t_i))}$, $w_i\in W_i \cap f^{-1}(\phi(t_i))$.
\end{itemize}
Since $\alpha$ is winning, there exists $w\in \bigcap_{i=1}^\infty W_i$. Let $\mathscr N$ be the set of pairs $(N,k)$, where $N$ is a neighborhood of
$w$ contained in $W_k$. For every such pair $N \cap f^{-1}(\phi(t_k)) \not=\emptyset$ by iii. We pick $w(N,k) \in N$ with $f(w(N,k))\in \phi(t_k)$. 
Consider $\mathscr N$ directed by the relation $(N,n) \preceq (N',n')$ iff $N' \subseteq N$ and $n' \geq n$. Then the net
$\mathcal N=\langle w(N,n): (N,n) \in \mathscr N\rangle$ converges to $w$, while the net $f(\mathcal N)$ converges to a point $v\in \bigcap_{i=1}^\infty \phi(t_i)$.
Since the graph of $f$ is closed, we deduce $v=f(w)$. But $f(w) \in \gamma(r_i) \subset G_i$ for every $i$  implies $f(w)\in G$,
while $f(w) \in \phi(t_1)\subset V$ gives $f(w) \in V$. On the other hand, $w\in W_k \subset H_k$ for every $k$ gives $w\in H$. We have shown $w\in W \cap f^{-1}(V\cap G)\cap H$.
\hfill $\square$
\end{proof}


\begin{theorem}
\label{theorem8}
Let $E$ be Baire, $F$ regular and $c$-complete, $f:E \to F$ a nearly continuous and nearly open mapping with closed graph.
Then $f$ is continuous.
\end{theorem}

\begin{proof}
By assumption $F$ admits a web $(T,\phi)$ satisfying property $(c)$ in the definition of $c$-completeness; cf. Section \ref{prepare}.
Now let $x\in E$ and $U$ a neighborhood of $f(x)$. By regularity of $F$ choose an open neighborhood $V$ of $f(x)$
with $f(x) \in V \subset \xoverline[.9]{V} \subset U$. Since $\xoverline[.9]{f^{-1}(V)}$ is a neighborhood of $x$, it suffices to prove
$\xoverline[.9]{f^{-1}(V)} \subset f^{-1}(U)$. Let $y\in \xoverline[.9]{f^{-1}(V)}$, it suffices to prove $f(y)\in \xoverline[.9]{V}$. 
To this end  let $W$ be an open neighborhood
of $f(y)$. It remains to prove $V\cap W \not=\emptyset$.
We define a strategy $\beta$ in the BM-game on $E$.

We have to define $\beta(\emptyset)$. We have
$\xoverline[.9]{f^{-1}(V)} = \xoverline[1.0]{f^{-1}(\bigcup \{ \phi(t): \phi(t) \subset V\})}$ by Lemma \ref{help}, and since $\xoverline[.9]{f^{-1}(W)}^\circ$ is a neighborhood
of $y$, we have $\xoverline[.9]{f^{-1}(W)}^\circ\cap f^{-1}(\bigcup \{ \phi(t): \phi(t) \subset V\}) \not=\emptyset$.
We choose $t_1\in T$ and $z_1 \in \xoverline[.9]{f^{-1}(W)}^\circ$ with $f(z_1) \in \phi(t_1) \subset V$. Now
$\xoverline[.9]{f^{-1}(W)} = \xoverline[1.0]{f^{-1}(\bigcup\{\phi(s): s\in T, \phi(s) \subset W\})}$ by Lemma \ref{help}, and since
$\xoverline[.9]{f^{-1}(\phi(t_1))}^\circ$ is a neighborhood of $z_1$,  it intersects
$f^{-1}(\bigcup\{\phi(s): s\in T, \phi(s) \subset W\})$. We choose $s_1 \in T$ with $\phi(s_1)\subset W$ and $y_1 \in \xoverline[.9]{f^{-1}(\phi(t_1))}^\circ$
such that $f(y_1)\in \phi(s_1) \subset W$. Then the open $U_1 =  \xoverline[.9]{f^{-1}(\phi(t_1))}^\circ \cap \xoverline[.9]{f^{-1}(\phi(s_1))}^\circ$ is nonempty,
and we let $\beta(\emptyset)=U_1$.

Now let $U_1' \subset U_1$ be nonempty open. We have to define $\beta(U_1,U_1')$. Since $U_1' \subset U_1 \subset \xoverline[.9]{f^{-1}(\phi(t_1))}^\circ$
and 
$\xoverline[.9]{f^{-1}(\phi(t_1))} = \xoverline[1]{f^{-1}(\bigcup\{\phi(t): t_1 <_Tt\})}$ by Lemma \ref{help}, $U_1'$ intersects $f^{-1}(\bigcup\{\phi(t): t_1 <_Tt\})$.
Choose $t_2\in T$ with $t_1 <_Tt_2$ and  $z_2 \in U_1' \cap f^{-1}(\phi(t_2))$. Now $f(z_2) \in \phi(t_2)$, hence
$\xoverline[.9]{f^{-1}(\phi(t_2))}^\circ$ is a neighborhood of $z_2 \in U_1' \subset \xoverline[.9]{f^{-1}(\phi(s_1))}^\circ$. But
$\xoverline[.9]{f^{-1}(\phi(s_1))} = \xoverline[1]{f^{-1}(\bigcup\{\phi(s): s_1 <_T s\})}$ by Lemma \ref{help}, hence $U_1'\cap \xoverline[.9]{f^{-1}(\phi(t_2))}^\circ$ intersects
$f^{-1}(\bigcup\{\phi(s): s_1 <_Ts\})$. Find $s_2\in T$ with $s_1 <_Ts_2$ and $y_2 \in U_1'\cap \xoverline[.9]{f^{-1}(\phi(t_2))}^\circ$ such that $f(y_2)\in \phi(s_2)$.
It follows that $U_2 = U_1'\cap \xoverline[.9]{f^{-1}(\phi(t_2))}^\circ \cap \xoverline[.9]{f^{-1}(\phi(s_2))}^\circ$ is nonempty, and we put
$\beta(U_1,U_1') = U_2$.

Continuing to define $\beta$ in this way, let $\alpha$ be a strategy winning against $\beta$. Then their play $U_1 \supset U_1' \supset U_2 \supset U_2' \supset \dots$
has the following properties:
\begin{enumerate}
\item $U_k' \subset U_k=U_{k-1}' \cap \xoverline[.9]{f^{-1}(\phi(t_k))}^\circ \cap \xoverline[.9]{f^{-1}(\phi(s_k))}^\circ$.
\item $z_{k+1} \in U_k' \cap f^{-1}(\phi(t_{k+1}))$, $t_k <_Tt_{k+1}$.
\item $y_{k+1} \in U_k' \cap \xoverline[.9]{f^{-1}(\phi(t_{k+1}))}^\circ$, $f(y_{k+1}) \in \phi(s_{k+1})$, $s_k <_Ts_{k+1}$.
\end{enumerate}
Since $\alpha$ beats $\beta$, there exists
$u\in \bigcap_{i=1}^\infty U_i \not=\emptyset$. Let $\mathscr N$ be the set of all pairs $(N,k)$ where $N$
is a neighborhood
of $u$ contained in $U_k$, ordered by the relation $(N,k) \preceq (N',k')$ iff $N'\subseteq N$ and $k'\geq k$. For $(N,k)\in \mathscr N$ we have $N \cap f^{-1}(\phi(t_k)) \not=\emptyset$, and also $N \cap f^{-1}(\phi(s_k))\not=\emptyset$, due to property 1. above. Pick $v(N,k) \in N \cap f^{-1}(\phi(t_k))$ and
$w(N,k) \in N \cap f^{-1}(\phi(s_k))$. Clearly the nets $\langle v(N,k)\rangle$ and $\langle w(N,k)\rangle$ both converge to $u$. 
On the other hand, by the completeness property $(c)$ of the web $(T,\phi)$ the net $\langle f(v(N,k)) \rangle$ has a cluster point
$v\in \bigcap_{i=1}^\infty \phi(t_i) \subset V$ and similarly, the net $\langle f(w(N,k))\rangle$ has a cluster $w\in \bigcap_{i=1}^\infty \phi(s_i) \subset W$,
so $v\in V$ and $w\in W$. But by closedness of the graph we have $(u,v)\in {\rm graph}(f)$ and $(u,w) \in {\rm graph}(f)$, hence $v=w$,
proving $V \cap W \not=\emptyset$.
\hfill $\square$
\end{proof}



\end{document}